\documentclass{article}
\usepackage{booktabs}
\usepackage{titlesec}

\titleformat{\section}
  {\fontsize{13pt}{15pt}\selectfont\bfseries}{\thesection}{1em}{}  
\titleformat{\subsection}
  {\fontsize{12pt}{15pt}\selectfont\bfseries}{\thesubsection}{1em}{}

\usepackage{multirow}
\usepackage{booktabs}
\usepackage{setspace}
\usepackage{mathrsfs}
\usepackage{todonotes}
\usepackage{geometry}
\usepackage{amsmath}
\usepackage{amssymb}
\usepackage{amsthm}
\usepackage{bm}
\usepackage{graphicx}
\graphicspath{{./Figure/}}
\usepackage{color}
\usepackage{tabu}
\usepackage{booktabs}
\usepackage{caption}
\usepackage{authblk}
\usepackage{subfigure}
\usepackage{multirow}

\newtheorem{remark}{Remark}
\newtheorem{definition}{Definition}

\begin{document}
\title{AP-MIONet: Asymptotic-preserving multiple-input neural operators for capturing the high-field limits of collisional kinetic equations}

\author[1 \thanks{Corresponding author: zhangtianai@sjtu.edu.cn}]{Tian-ai Zhang}
\author[1, 2, 3]{Shi Jin}
\affil[1]{School of Mathematical Sciences, Shanghai Jiao Tong University, Shanghai, China}
\affil[2]{Institute of Natural Sciences, MOE-LSC, Shanghai Jiao Tong University, Shanghai, China}
\affil[3]{Shanghai Artificial Intelligence Laboratory, Shanghai, China}
\date{\today}
\maketitle

\begin{abstract}
In kinetic equations, external fields play a significant role, particularly when their strength is sufficient to balance collision effects, leading to the so-called high-field regime. Two typical examples are the Vlasov-Poisson-Fokker-Planck (VPFP) system in plasma physics and the Boltzmann equation in semiconductor physics. In this paper, we propose a generic asymptotic-preserving multiple-input DeepONet (AP-MIONet) method for solving these two kinetic equations with variable parameters in the high-field regime. Our method aims to tackle two major challenges in this regime: the additional variable parameters introduced by electric fields, and the absence of an explicit local equilibrium, which is a key component of asymptotic-preserving (AP) schemes. We leverage the multiple-input DeepONet (MIONet) architecture to accommodate additional parameters, and formulate the AP loss function by incorporating both the mass conservation law and the original kinetic system. This strategy can avoid reliance on the explicit local equilibrium, preserve the mass and adapt to non-equilibrium states. We demonstrate the effectiveness and efficiency of the proposed method through extensive numerical examples.

\vspace{1em}\noindent\textbf{Keywords:} Multiscale kinetic equations, High-field limit, Asymptotic-preserving schemes, Physics-informed multiple-input neural operators

\noindent\textbf{Mathematics Subject Classification:} 82D10, 82B40, 82C40, 68T07
\end{abstract}

\section{Introduction}
Kinetic equations are fundamental multiscale models bridging continuum and atomistic models in various physical realms, particularly in plasma and semiconductor physics \cite{liboff2003kinetic, jin2022asymptotic}. These equations characterize the statistical behavior of particle systems by modeling the transport and collisions of particles, as well as their interactions with external fields or surrounding media. In the presence of external fields, the high-field regime is a common scenario where the external field is strong to balance the collision term \cite{nieto2001high}. The asymptotic behaviors in this regime, referred to as the high-field limits, have been widely studied in the fields of plasma and semiconductor physics \cite{nieto2001high, goudon2005multidimensional, frosali1989scattering, frosali1989conditions, poupaud1992runaway, cercignani1997high, abdallah2001high}. The Vlasov-Poisson-Fokker-Planck (VPFP) system and semiconductor Boltzmann equation are classical kinetic models for describing the high-field limit in plasma and semiconductor physics, respectively. Modelling and simulating electron dynamics in plasma and semiconductor devices at the kinetic level pose significant challenges due to the high dimensionality of the phase space (typically six dimensions plus time) and the presence of multiple spatial and temporal scales.

In recent years, neural network (NN) methods have developed to pivotal techniques for solving partial differential equations (PDEs), leveraging the powerful approximation capabilities of NNs \cite{hornik1990universal}. Extensive research has also been conducted on multiscale kinetic equations and hyperbolic systems \cite{chen2022solving, hwang2020trend, xiao2023relaxnet, jin2023asymptotic1, liu2020multi,lou2021physics}. Deep learning techniques show great potential in solving high-dimensional problems \cite{han2018solving}, an inherent challenge for classical numerical methods. Their mesh-free nature facilitates the solution of PDEs defined on complex domains and geometries with ease of implementation. While conventional NN architectures are proficient in learning mappings between finite-dimensional Euclidean spaces, they exhibit limited generalization capabilities when facing with equations involving multiple scales and other physical constraints, structures and parameters.
In particular,  solving these so-called \textit{parametric PDE} problems is time-consuming, since we need to re-train the NNs whenever parameters of a given PDE problem are changed.

To enhance computational efficiency, neural operators \cite{goswami2023physics} have been proposed to directly learn the \textit{solution operator} that maps variable input functions to their corresponding solutions of the governing PDE system across infinite-dimensional Banach spaces. These operators allow the neural network to be trained \textit{only once} to generate solutions for a class of PDEs. The Deep Operator Network (DeepONet) \cite{lu2021learning} is a pioneering work in the field of neural operators, grounded in rigorous approximation theory. Numerous neural operators have been proposed, demonstrating superior performance across a diverse range of applications \cite{li2020fourier, li2024physics, li2020neural, you2022nonlocal, li2020multipole, jin2022mionet, fanaskov2023spectral}. In the domain of kinetic equations, neural operators have showed great promise for building fast emulators to solve parametric PDE systems \cite{lee2024structure, lee2023oppinn, chen2022solving, xu2023transfer, wu2024capturing}.

Applying neural operators to solve kinetic equations with the high-field scaling presents two significant challenges: managing multiple input functions for solution operators, and tackling inherent multiscale nature of kinetic equations. For the first challenge, neural operators are typically designed to learn operators defined within a single Banach space, thereby limiting them to processing only one input function. However, the presence of an external field, governed by the coupled equation, could involve additional varying parameters as inputs for the solution operators besides IBCs. The multiple-input DeepONet (MIONet) \cite{jin2022mionet} was proposed to map from a product of multiple Banach spaces to another Banach space, enabling it to accept multiple input functions. This model generalizes the DeepONet in both theory and architecture, making it possible to learn a broad range of more practical and complex operators \cite{jiang2023fourier, wang2024spi, kong2023b, hu2024hybrid}.

For the second multiscale challenge, vanilla NN and neural operator methods often fail to capture the correct asymptotic limit behavior of kinetic equations when scale parameters are small \cite{ jin2023asymptotic, wu2024capturing}. The asymptotic-preserving (AP) mechanism, a classical computational paradigm for addressing multiscale problems, has been successfully adapted to the domains of NNs and neural operators \cite{li2022apfos, jin2023asymptotic, lu2022solving, bertaglia2022asymptotic, li2022model, jin2023asymptotic1, wu2024capturing}. An AP computational method can preserve the continuous asymptotic limits along the transition from microscopic to macroscopic models in a numerically stable way \cite{jin2022asymptotic}. This scheme allows efficient numerical approximations across all regimes without numerically resolving small scales. Within the field of operator learning, Wu \textit{et al.} \cite{wu2024capturing} first introduced the concept of an AP loss function for the linear transport equation with diffusive scaling, leading to the so-called Asymptotic-Preserving Convolutional DeepONets (APCONs). This approach reformulated the vanilla loss function, such that the loss of the microscopic equation transits to the loss of the corresponding macroscopic equation at the asymptotic limit.

This study aims to apply neural operator methods to solve multiscale kinetic equations in the high-field regime. Our focus is on two classical equations: the VPFP system and the semiconductor Boltzmann equation. To accommodate multiple input functions derived from the external field and IBCs, we utilize MIONets to parameterize the solution operators of underlying PDE systems. To takle the multiscale challenge, we incorporate the AP mechanism into the loss function of the MIONets. Specifically, we modify the loss function within the physics-informed MIONets (PI-MIONets) \cite{zheng2023state} framework, where the vanilla loss is constructed by taking $L^2$-residuals of underlying PDE constraints. A remarkable difficulty arises from the frequent absence of an explicit local equilibrium in the high-field regime,  particularly for the semiconductor Boltzmann equation. This difficulty makes the modern AP loss designs \cite{jin2023asymptotic1, jin2023asymptotic2, wu2024capturing, lu2022solving} intractable to implement, which are usually based on the micro-macro and even-odd decompositions. Our key idea is to incorporate both the mass conservation law and the original kinetic equations into the loss function. This strategy also ensures that the loss preserves the mass and fits for non-equilibrium situations, both of which are crucial but challenging for neural networks \cite{jin2023asymptotic2, jin2023asymptotic}. As a result, we propose a new Asymptotic-preserving multiple-input DeepONet (AP-MIONet) generic to both the VPFP system and the semiconductor Boltzmann equation with the high-field scaling. This method extends the AP concept within the operator learning framework for a practical and complex scenario, where external fields play an important role in kenetic models.

This paper is organized as follows: Sect. 2 introduces the VPFP system in plasma physics, the semiconductor Boltzmann equation in semiconductor physics and their corresponding high-field limits. In Sect. 3, we review the PI-MIONets equipped with modified multi-layer perceptrons (MLPs), whose loss functions are limited in accurately capturing the high-field limit equations. In Sect. 4, we propose the new AP-MIONet method and formally validate the AP property of its loss function. Numerous numerical examples are presented in Sect. 5 to demonstrate the effectiveness and efficiency of the AP-MIONet method. We conclude the paper in Sect. 6.

\section{Introduction to two kinetic equations and their high-field limits}\

Consider the kinetic system that  describe the time evolution of electron density distribution function $f(t, x, v)\geq 0$ under the action of a self-consistent potential $\phi$: 
\begin{subequations}
\begin{equation}\label{transport}
    \partial_t f+v \cdot \nabla_x f-\frac{1}{\varepsilon} \nabla_x \phi \cdot \nabla_v f =\frac{1}{\varepsilon} \mathscr{Q}(f),
\end{equation}
\begin{equation}\label{poisson}
    -\triangle_x \phi(t, x)=\rho(t, x)-h(x).
\end{equation}
\end{subequations}
Here, the spatial variable $x$ and velocity $v$ belong to $\mathbb{R}^N$, where $N$ is the dimension of the system. The time $t$ belongs to $\mathbb{R}_{+}$. The parameter $\varepsilon$ characterizes the scaled mean free path of electrons. The electric potential $\phi$ is governed by the Poisson equation \eqref{poisson}. The electron density $\rho$ is defined as 
\[
\rho(t, x)=\int_{\mathbb{R}^N} f(t, x, v) \mathrm{d} v.
\]
The function $h(x)$
represents a given positive background charge density. The global neutrality of the systems can be satisfied by the following relation:
\begin{equation}
    \int_{\mathbb{R}^N} \int_{\mathbb{R}^N} f(t=0, x, v) \mathrm{d} x \mathrm{~d} v=\int_{\mathbb{R}^N} h(x) \mathrm{d} x .
\end{equation}
Different collision operators $\mathscr{Q}$ are applied to characterize various physical systems. In plasma physics, the Fokker-Planck operator is adopted, leading to the Vlasov-Poisson-Fokker-Planck (VPFP) system \cite{degend1986global}. The semiconductor Boltzmann operator, regarded as a more general collision operator, is used to describe the semiconductor kinetic model \cite{poupaud1992runaway, abdallah2001high}.
Under the high-field scaling, the effect of electric field is strong  compared with collisions. The asymptotic procedure of this case is called the high-field limit.  
A detailed introduction to these two equations and their high-field limits is provided in the following sections.

\subsection{The VPFP system and its high-field limit}

The VPFP system is a kinetic description of the Brownian motion for a large system of particles in a surrounding bath \cite{chandrasekhar1943stochastic}. In electrostatic plasma, it models the interactions between electrons and the surrounding bath through the Coulomb force. Specifically, this system governs the time evolution of the electron distribution function $f:(t, x, v) \in \mathbb{R}_{+} \times \mathbb{R}^N \times \mathbb{R}^N \rightarrow \mathbb{R}_{+}$ acted by the potential $\phi(t,x)$:
\begin{subequations}
    \begin{equation}\label{vlasov}
        \partial_t f+v\cdot\nabla_x f-\frac{1}{\varepsilon}\nabla_x \phi \cdot \nabla_v f =\frac{1}{\varepsilon} \nabla_{v} \cdot\left[v f+\nabla_{v} f\right] := \frac{1}{\varepsilon} \mathscr{Q}_\text{FP}(f),
    \end{equation}
    \begin{equation}\tag{1b}
        -\triangle_x \phi(t, x) =\rho(t, x)-h(x).
    \end{equation}
\end{subequations} 
Here, the operator $\mathscr{Q}_\text{FP}$ represents a linear Fokker-Planck operator, defined as:
\begin{equation}\label{fpop}
    \mathscr{Q}_\text{FP}(f)(t, x, v):=\nabla_v\cdot\left[v f(t, x, v)+\nabla_v f(t, x, v)\right].
\end{equation}
The parameter $\varepsilon:=\left(l_e/\Lambda\right)^2$ is the ratio between the mean free path of electrons $l_e$ and the Debye length $\Lambda$ \cite{livadiotis2014electrostatic}. Under this scaling, the electric field is strong and comparable to the collision. The limiting process of $\varepsilon \rightarrow 0$ leads to the high-field limit of the VPFP system.

We now formally derive the limit equation of the VPFP system. First, integrating the Vlasov equation \eqref{vlasov} over $v$ in $\mathbb{R}^N$ gives:
\begin{equation}\label{deriv1}
    \partial_t \int_{\mathbb{R}^N} f d v+\nabla_{x} \cdot \int_{\mathbb{R}^N} v f d v-\frac{1}{\varepsilon} \int_{\mathbb{R}^N} \nabla_v \cdot\left(\nabla_{x} \phi f d v\right)=\frac{1}{\varepsilon} \int_{\mathbb{R}^N} \nabla_{v} \cdot\left(v f+\nabla_{v} f\right) d v.
\end{equation}
After integrating by parts, the continuity equation can be derived:
\begin{equation}\label{rhotrans}
    \partial_t \rho+\nabla_x \cdot j=0,
\end{equation}
where the flux $j$ is defined as 
\[
j:=\int_{\mathbb{R}^N} v f(t, x, v) \mathrm{d} v.
\]
Next, one multiplies the Vlasov equation \eqref{vlasov} by $v$, integrates over $v$ in $\mathbb{R}^N$ and takes the limit as $\varepsilon \rightarrow 0$, yielding:
\begin{equation}\label{deriv2}
0=\int_{\mathbb{R}^N} f \nabla_{x} \phi+v f+\nabla_v f \mathrm{d} v.
\end{equation} 
Therefore, one has:
\begin{equation}\label{deriv3}
    j=-\rho\left(\nabla_x \phi\right).
\end{equation}
Substituting Eq.\eqref{deriv3} in Eq.\eqref{rhotrans} leads to the high-field limit system:
\begin{equation}\label{limit}
    \left\{\begin{array}{l}
    \partial_t \rho-\nabla_x \cdot\left(\rho \nabla_x \phi\right)=0,\\
    -\triangle_x \phi=\rho-h(x).
    \end{array}\right.
\end{equation}
A rigorous proof for the high-field limit in the one-dimensional case was presented in \cite{nieto2001high}. Goudon \textit{et al.} \cite{goudon2005multidimensional} extended this proof to the multi-dimensional case for the electrostatic VPFP system.

A key feature of the VPFP system is its explicit expression for the local equilibrium, a critical component to formulate asymptotic-preserving (AP) methods \cite{crouseilles2011asymptotic, jin2011asymptotic}. Here, we derive this explicit equilibrium. By combining the force term $\frac{1}{\varepsilon} \nabla_x \phi \cdot \nabla_v f$ and the collision term $\frac{1}{\varepsilon}\mathcal{Q}(f)$ in the Vlasov equation \eqref{vlasov}, the VPFP system \eqref{vlasov}-\eqref{poisson} can be rewritten as:
\begin{equation}\label{newvfp}
    \begin{aligned}
        \partial_t f+v \cdot \nabla_x f &=\frac{1}{\varepsilon} \nabla_v\cdot\left[(v+\nabla_x\phi) f+\nabla_v f\right] := \frac{1}{\varepsilon} \mathscr{L} f,\\
        -\Delta_x \phi &=\rho-h(x).
    \end{aligned}
\end{equation}
Here, $\mathscr{L}$ is a linear operator that depends on $\nabla_x \phi$ and retains the properties of the Fokker-Planck operator. The null space and the rank of the operator $\mathscr{L}$ are 
\[
\mathscr{N}(\mathscr{L})=\operatorname{Span}\{\mathscr{M}\}=\{f=\rho \mathscr{M}\},
\]
and 
\[
\mathscr{R}(\mathscr{L})=(\mathscr{N})^{\perp}(\mathscr{L})=\{f\, \text{such that} \, \langle f\rangle:=\int_{\mathbb{R}^N} f d v=0\}
\]
respectively. Here, $\mathscr{M}$ is the so-called \textit{local Maxwellian}, depending on $(t,x)$ through the potential $\phi$:
\begin{equation}\label{vlasov2}
    \mathscr{M}(t, x, v)=\frac{1}{\sqrt{2 \pi}} \exp\left(-\frac{|v+\nabla_x\phi(t, x)|^2}{2}\right).
\end{equation}
It can be formally verified that as $\varepsilon\rightarrow0$ in Eqs.\eqref{newvfp}, $f$ converges to $\rho\mathscr{M}$ with $\rho$ satisfying the high-field limit equation \eqref{limit}.

\subsection{The semiconductor Boltzmann equation and its high-field limit}

In the semiconductor kinetic model, the Boltzmann equation (BE) describes the carrier transport in semicondutor devices.  Consider  the semi-classical evolution of the electron distribution function $f(t,x,v)$ under the parabolic band approximation is governed by the scaled Boltzmann equation:
\begin{equation}\label{semiB}
    \partial_t f+v \cdot \nabla_x f+\frac{1}{\varepsilon} E \cdot \nabla_v f=\frac{1}{\varepsilon}\mathscr{Q}_\text{semiB}(f), \quad t>0, x \in \mathbb{R}^N, v \in \mathbb{R}^N.
\end{equation}
 The electric field $E=-\nabla_x\phi$ is assumed to be either given or self-consistently determined through the Poisson equation \eqref{poisson}. The integral operator $\mathscr{Q}_\text{semiB}$ characterizes the interactions of the electrons with the semiconductor crystal lattice and between electrons themselves. For semiconductors with low electron density, the general form of $\mathscr{Q}$ is \cite{markowich2012semiconductor}:
\begin{equation}\label{nond1}
    \mathscr{Q}_\text{nond}(f)=\int_{\mathbb{R}^N}\left(s\left(v^{\prime}, v\right) f\left(t, x, v^{\prime}\right)-s\left(v, v^{\prime}\right) f(t, x, v)\right) d v^{\prime},
\end{equation}
which pertains to the non-degenerate semiconductor case. The cross section $s$ is a known positive function. In scenarios of high electron density in semiconductors, Pauli's exclusion principle becomes necessary. Hence, the collision operator $\mathscr{Q}_\text{semiB}$ is modified for the degenerate semiconductor case as follows:
\begin{equation}\label{deg1}
    \mathscr{Q}_\text{deg}(f)=\int_{\mathbb{R}^{N_v}}\left(s\left(v^{\prime}, v\right) f^{\prime}(1-f)-s\left(v, v^{\prime}\right) f\left(1-f^{\prime}\right)\right) d v^{\prime},
\end{equation}
where $f$ and $f^{\prime}$ are shorthanded notations for $f(t, x, v)$ and $f\left(t, x, v^{\prime}\right)$ respectively. The $(1-f)$ terms in Eq.\eqref{deg1} account for the Pauli exclusion principle, as well as introduce quadratic nonlinearity. The parameter $\varepsilon$ is the ratio of the mean free path to the typical length scale of semiconductor devices.

Numerically solving the BE is impractically prohibitive due to the high dimensionality of the equation, defined in the  phase space. To simplify computations, various macroscopic models based on the diffusion approximation have been developed, primarily including the Drift-Diffusion (DD) model \cite{poupaud1991diffusion, golse1992limite}, the Energy-Transport (ET) model \cite{ben1996energy, abdallah2000convergence} and the Spherical Harmonic Expansion (SHE) model \cite{stratton1957influence, stratton1962diffusion}. These macroscopic models are derived under conditions where collisions are at the dominant order. However, as semiconductor devices continue to miniaturize,  standard macroscopic models face challenges in accurately describing hot electron transport under strong electric fields ans one has to use kinetic models.  

Under the scaling in Eq.\eqref{semiB}, strong electric field can balance the scattering effects in the operator $\mathscr{Q}_\text{semiB}$. This scenario is commonly referred to as the high-field regime in semiconductor physics. The macroscopic limit, as $\varepsilon \rightarrow 0$, is termed as the high-field limit. This limit was first studied for the simplified linear Boltzmann equation \cite{frosali1989scattering, frosali1989conditions}, and later extended to the non-degenerate case \cite{poupaud1992runaway} and the degenerate case \cite{abdallah2001high}. When the electrostatic potential is self-consistently governed by the Poisson equation, the high-field limit was derived for the BGK-type collision \cite{cercignani1997high}. For more general collision operators, this topic remains an active area of research. In the following sections, we formally derive the high-field limits for both the non-degenerate and degenerate cases.

\subsubsection{The non-degenerate case}

We first simplify the operator $\mathscr{Q}_\text{nond}$ in the non-degenerate case. The cross section $s$ in $\mathscr{Q}_\text{nond}$ satisfies the principle of detailed balance \cite{blakemore2002semiconductor} yielding:
\begin{equation}\label{balance}
    s\left(v^{\prime}, v\right) M\left(v^{\prime}\right)=s\left(v, v^{\prime}\right) M(v),
\end{equation}
where the Maxwellian distribution $M(v)$ is defined as $M(v)=1/(2 \pi)^{N/2}\exp\left(-|v|^2/2\right)$. The null space of $\mathscr{Q}_\text{nond}$ is  spanned by $M(v)$. Considering the relation \eqref{balance}, it is advantageous to introduce a new function $\Psi$ defined as:
\begin{equation}
    \Psi\left(v, v^{\prime}\right)=\frac{s\left(v^{\prime}, v\right)}{M(v)}, \quad \text { so that } \Psi\left(v, v^{\prime}\right)=\Psi\left(v^{\prime}, v\right).
\end{equation}
This $\Psi$ enables the operator $\mathscr{Q}_\text{nond}$ to be reformulated in a symmetric form:
\begin{equation}\label{nond2}
    \mathscr{Q}_\text{nond}(f)=\int_{\mathbb{R}^N} \Psi\left(v, v^{\prime}\right)\left(M(v) f\left(t, x, v^{\prime}\right)-M\left(v^{\prime}\right) f(t, x, v)\right) d v^{\prime} .
\end{equation}
For details, readers are directed to \cite{blakemore2002semiconductor, markowich2012semiconductor}.

We now derive the high-field limit equation for this case. As $\varepsilon\rightarrow 0$ in Eq.\eqref{semiB}, the distribution $f$ converges to:
\begin{equation}\label{deriv4}
    f(t, x, v) \rightarrow \rho(t, x) F_{E(t, x)}(v),
\end{equation}
where $F_{E(t, x)}(v)$ is the solution to:
\begin{equation}\label{nond_eq2}
    \int_{\mathbb{R}^N} F_E(v) d v=1, \quad E \cdot \nabla_v F_E-\mathscr{Q}_\text{nond}\left(F_E\right)=0, \quad F_E \geq 0.
\end{equation}
Integrating Eq.\eqref{semiB} over $v$ yields the equation for the macroscopic density $\rho$: 
\begin{equation}\label{deriv5}
    \partial_t \rho(t, x)+\int_{\mathbb{R}^{d_v}} v \cdot \nabla_x f=0.
\end{equation}
Substituting Eq.\eqref{deriv4} to Eq.\eqref{deriv5}, one derives the corresponding limit equation:
\begin{equation}\label{nond_eq1}
    \partial_t \rho(t, x)+\nabla_x \cdot(\rho(t, x) \sigma(E(t, x)))=0, \quad \sigma(E)=\int_{\mathbb{R}^N} v F_E(v) d v.
\end{equation}
This limiting equation manifests as a linear convection equation for the macroscopic density, with the convection proportional to the scaled electric field. Eq.\eqref{nond_eq1} excludes the diffusion effect and essentially follows Ohm's law. The existence of a limit solution can be ensured by a criterion for the cross-section $s$ \cite{poupaud1992runaway}. It is important to note that not every $\mathscr{Q}_\text{nond}$ guarantees a unique solution for Eq. \eqref{nond_eq2}.

\subsubsection{The degenerate case}

Since the cross section $s$ in Eq.\eqref{deg1} satisfies the detailed balance like the non-degenerate case, the operator $\mathscr{Q}_\text{deg}$ can also be reformulated in a symmetric form:
\begin{equation}\label{deg}
    \mathscr{Q}_\text{deg}(f)(t, x, v)=\int_{\mathbb{R}^N} \Psi\left(v^{\prime}, v\right)\left(M(v) f\left(t, x, v^{\prime}\right)(1-f(t, x, v))-M\left(v^{\prime}\right) f(t, x, v)\left(1-f\left(t, x, v^{\prime}\right)\right)\right) d v^{\prime}.
\end{equation}
The null space of $\mathscr{Q}_\text{deg}(f)(t, x, v)$ is spanned by the Fermi-Dirac distribution:
\begin{equation}
    M_{F D}=\frac{1}{1+e^{\frac{m_e v^2}{2 K_B T}-\frac{\mu}{K_B T}}},
\end{equation}
where $T$ is temperature of the lattice, $m_e$ is the effective mass of electrons, $\mu$ is the electron Fermi energy and $K_B$ is the Boltzmann constant. In this case, the quantum effects are taken into account, resulting in the nonlinearity of the collision operator $\mathscr{Q}_\text{deg}$. Therefore, as $\varepsilon\rightarrow0$ in Eq.\eqref{semiB}, the distribution function $f$ can no longer be factorized into a product of two separate functions, where one function solely depends on the spatial and temporal variables $x$ and $t$, and the other function only depends on the velocity variable $v$.

We now derive the limit equation for the degenerate case. To ensure the existence of the limit solution, one considers the case where $\Psi$ belongs to $W^{2, \infty}\left(B^2\right)$ and satisfies $\Psi_0 \leq \Psi\left(v, v^{\prime}\right) \leq \Psi_1$ for some positive constants $\Psi_0$ and $\Psi_1$. Here, $B$ represents either the Brillouin zone or the entire space $\mathbb{R}^N$. As $\varepsilon$ tends to zero in Eq.\eqref{semiB}, the distribution function converges to the limit given by:
\begin{equation}\label{deriv6}
    f \rightarrow F(\rho(t, x), E(t, x))(v).
\end{equation}
Here, $F(\rho, E)(v)$ is the unique solution in space $D_E=\left\{F \in L^1(B) ; E \cdot \nabla_v F \in L^1(B)\right\}$, satisfying the following conditions:
\begin{equation}\label{cond7}
    E \cdot \nabla_v F-\mathscr{Q}_{deg}(F)=0, \quad \int_{\mathbb{R}^{d v}} F(t, x, v) d v=\rho(t, x), \quad 0 \leq F \leq 1.
\end{equation}
The mapping $(\rho, E) \mapsto F(\rho, E)$ from $\mathbb{R}^{+} \times \mathbb{R}^N$ to $L^1(B)$ is $C^2$ differentiable. Integrating Eq.\eqref{semiB} over $v$, one obtains the equation for the density $\rho$:
\begin{equation}\label{deriv7}
    \partial_t \rho(t, x)+\int_{\mathbb{R}^{d_v}} v \cdot \nabla_x f=0.
\end{equation}
Substituting Eq.\eqref{deriv6} to Eq.\eqref{deriv7}, the limit equation is derived:
\begin{equation}\label{deg_limit}
    \partial_t \rho(t, x)+\nabla_x(j(\rho(t, x) ; E(t, x)))=0,
\end{equation}
where $j(\rho ; E)=\int_{\mathbb{R}^d v} v F(\rho, E)(v) d v$. This result was established for a given $E(x) \in \mathbb{R}^N$ over local time intervals to obtain a regular limit solution \cite{bouchut1993existence}.

\begin{remark}
In the former non-degenerate case, the limit equation \eqref {nond_eq1} is linear in $\rho$, thereby ensuring the existence of a unique global solution in time. However, in the degenerate case, due to the nonlinearity of the flux function in Eq.\eqref{deg_limit}, the existence and uniqueness of regular solutions were only proved locally in time and shocks may subsequently develop \cite{ben2007high}. 
\end{remark}

We emphasize that, unlike the VPFP system, the semiconductor Boltzmann equation does not possess an explicit expression for the local equilibrium in the high-field limit. This absence poses a challenge in the design of neural network-based asymptotic-preserving (AP) methods. Therefore, the primary goal of
this work is to develop an effective AP method within the operator learning framework without reliance on the explicit local equilibrium.

\section{Operator learning for the two kinetic equations}
This section provides a concise overview of physics-informed multiple-input DeepONets (PI-MIONets) equipped with modified multi-layer perceptrons (MLPs). We explore to apply PI-MIONets to solve the two kinetic equations in the high-field regime, where the electric potential is governed by the Poisson equation. We discuss the limitations of PI-MIONets in the case with a small scale parameter.

\subsection{A primer on PI-MIONets with modified MLPs}
\noindent $\bullet$ \textbf{Multiple-Input DeepONets} \hspace{0.5em} DeepONets are an emerging paradigm in learning solution operators that map between infinite-dimensional function spaces for a wide range of dynamic systems and PDEs. This approach is based on the universal approximation theorem for operators \cite{chen1995universal}, which ensures that a neural network with sufficient width can approximate any continuous function with the desired degree of precision. The architecture of a DeepONet consists of a branch net and a trunk net for encoding the input function and  the domain of the output function, respectively. However, this original architecture is designed to learn operators defined on a single Banach space, limiting its input to only one function. This limitation prevents DeepONets from supporting complex operators in realistic setups. For example, solution operators for PDE systems typically map from both the initial and boundary conditions, which often reside in two distinct domains. To accommodate more input functions, the multiple-input DeepONet (MIONet) \cite{jin2022mionet} was developed to learn operators defined on the tensor product of multiple Banach spaces. It generalizes the DeepONet both theoretically and numerically.

\begin{figure}[!htb]
	\centering  
	\subfigbottomskip=5pt
	\subfigcapskip=-5pt
	\subfigure[The architecture of a MIONet.]{
	\includegraphics[width=0.55\linewidth]{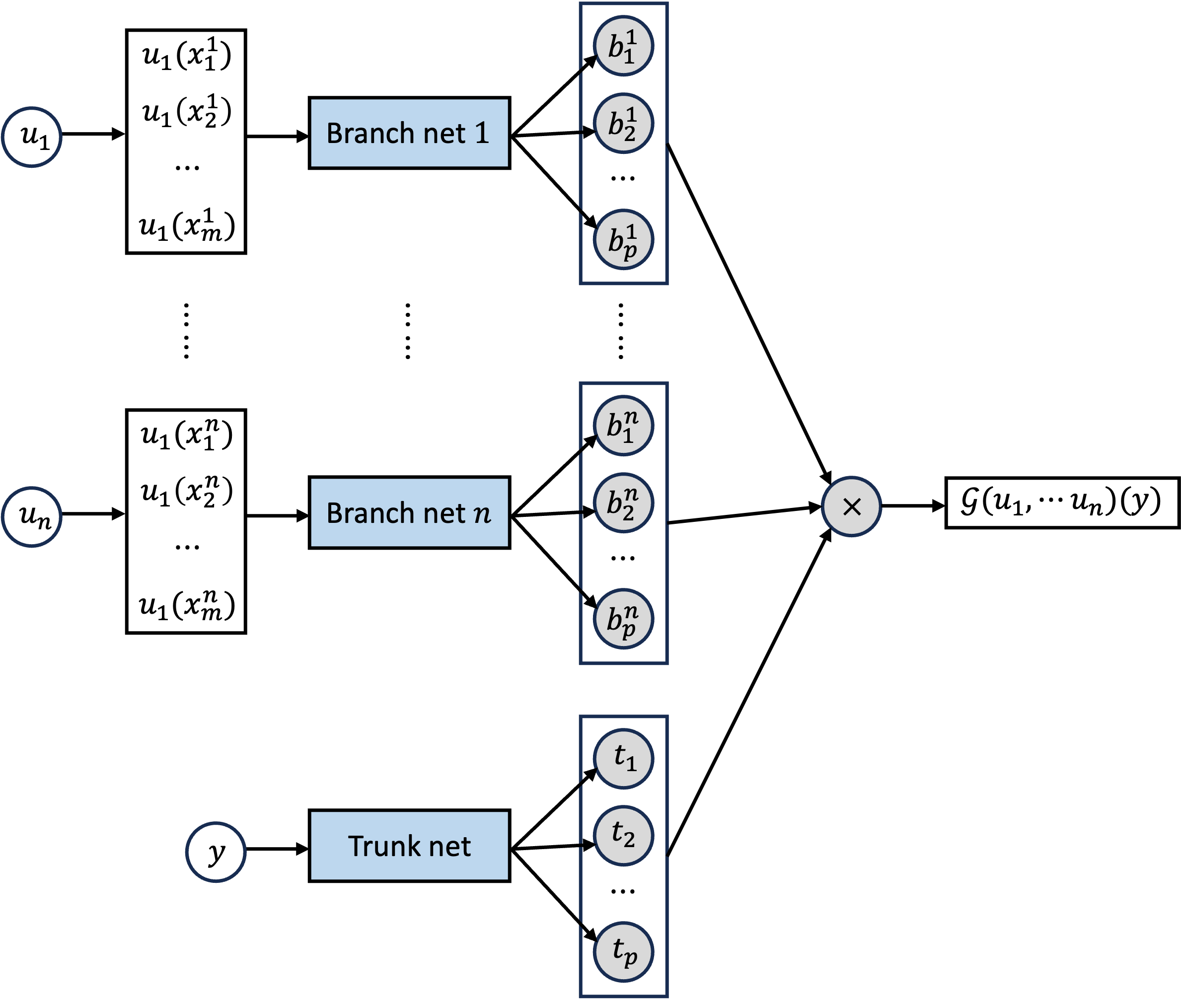}}
 
	\subfigure[The architecture of a PI-MIONet.]{
	\includegraphics[width=0.70\linewidth]{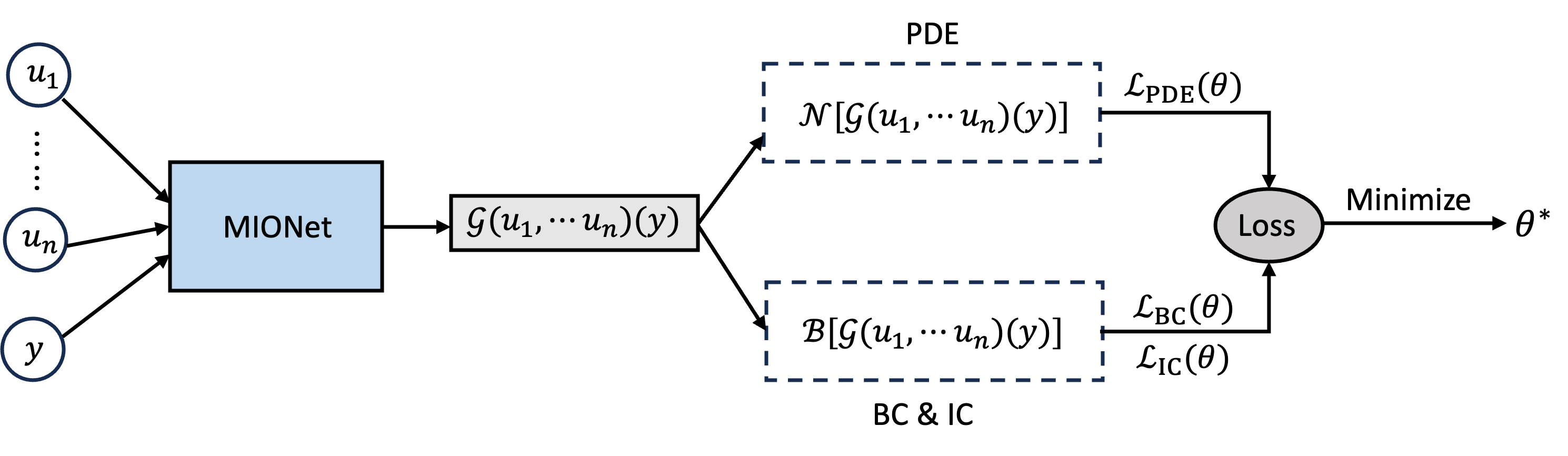}}
 
  \caption{The architecture of (a) MIONet and (b) physics-informed MIONet (PI-MIONet). In subfigure (a), assume the operator $\mathcal{G}$ accepts $n$ input functions $\boldsymbol{u}_1, \cdots \boldsymbol{u}_n$ and is evaluated at coordinate $\boldsymbol{y}$. The MIONet architecture corresponding to $\mathcal{G}$ comprises $n+1$ sub-networks: the $n$ branch nets, each for extracting latent representations from their respective input function $\{\boldsymbol{u}_i\}_{i=1}^n$; and a single trunk net for extracting latent representations from the input coordinate $\boldsymbol{y}$ at which the output function is evaluated. A continuous and differentiable representation of the output function $\mathcal{G}(\boldsymbol{u}_1,\cdots, \boldsymbol{u}_n)(\boldsymbol{y})$ is obtained by integrating the latent representations extracted by all sub-networks through a dot product. In subfigure (b), to make the MIONet physics-informed, additional regularization mechanisms are incorporated via automatic differentiation for biasing the MIONet output to satisfy the governing PDE system, boundary condition (BC) and initial condition (IC).}
\label{PI-MIONet1}
\end{figure}

To elaborate on the architecture of MIONets, let us consider an operator $\mathcal{G}$ that receives $n$ input functions $\boldsymbol{u}_i$ for $i=\{1,\cdots, n\}$, with $\mathcal{G}(\boldsymbol{u}_1,\cdots,\boldsymbol{u}_n)$ denoting the corresponding output function. The aim of MIONets is to approximate the operator $\mathcal{G}(\boldsymbol{u}_1,\cdots,\boldsymbol{u}_n)$ evaluated at continuous coordinates $\boldsymbol{y}$. The operator $\mathcal{G}$ is parameterized by a MIONet $\mathcal{G}_\theta$, where $\theta$ represents all trainable parameters in the MIONet \cite{jin2022mionet}. The architecture of this MIONet comprises $n$ independent branch nets and a single trunk net. The $i$-th branch net takes $\boldsymbol{u}_i$ as input and returns  features embedding $\left[b_1^i, b_2^i, \ldots, b_p^i\right]^T \in \mathbb{R}^p$ as output, where $\boldsymbol{u}_i=\left[\boldsymbol{u}_i\left(\boldsymbol{x}_1^i\right), \boldsymbol{u}_i\left(\boldsymbol{x}_2^i\right), \ldots, \boldsymbol{u}_i\left(\boldsymbol{x}_m^i\right)\right]$ denotes the function $\boldsymbol{u}_i$ evaluated at a set of fixed sensors $\left\{\boldsymbol{x}_k^i\right\}_{k=1}^m$ (where the input functions are defined). The trunk net inputs the continuous coordinates $\boldsymbol{y}$ for $\mathcal{G}_\theta(\boldsymbol{u}_1,\cdots,\boldsymbol{u}_n)$ evaluation and outputs  features embedding $\left[t_1, t_2, \ldots, t_p\right]^T \in \mathbb{R}^p$. The final output of the MIONet is computed as a dot product of the outputs from the branch and trunk nets:
\begin{equation}
\mathcal{G}_\theta\left(\boldsymbol{u}_1, \boldsymbol{u}_2, \ldots, \boldsymbol{u}_n\right)(\boldsymbol{y})=\sum_{j=1}^p \underbrace{b_j^1\left(\boldsymbol{u}_1\right)}_{\text {branch }1} \times \underbrace{b_j^2\left(\boldsymbol{u}_2\right)}_{\text {branch }2} \cdots \times \underbrace{b_j^n\left(\boldsymbol{u}_n\right)}_{\text {branch }n} \times \underbrace{t_j(\boldsymbol{y})}_{\text {trunk }}+b_0,
\end{equation}
where $b_0$ is a bias term added to enhance the performance of the MIONet. A diagram of this MIONet architecture is displayed in the top panel of Fig. \ref{PI-MIONet1}.

\vspace{1em}\noindent $\bullet$ \textbf{Physics-informed MIONets} 

\hspace{0.5em} For MIONet training, purely data-driven approaches face two challenges. First, it relies on an extensive dataset to generalize solutions, which increases the burden on data collection and requires substantial memory allocation. Compared with DeepONets, MIONets demand a larger volume of training samples since the learned operators are defined across multiple function spaces. The second challenge arises because solution operators trained by data-driven methods is only a rough approximation and may not accurately satisfy the underlying PDEs. To address these challenges, the physics-informed MIONet (PI-MIONet) \cite{zheng2023state} was developed. PI-MIONet leverage automatic differentiation \cite{paszke2017automatic} to integrate physical governing laws in its loss function by penalizing residuals of PDEs, since the output of the MIONet is differentiable with respect to the input coordinates in its trunk net. By minimizing the physics-informed loss, the target output function is forced to align with underlying physical constraints using only appropriate initial and boundary conditions, thereby eliminating the need for paired input-output data \cite{wang2021learning}. A schematic diagram of the PI-MIONet architecture is depicted in the bottom panel of Fig. \ref{PI-MIONet1}.

\vspace{1em}\noindent $\bullet$ \textbf{Modified multi-layer perceptrons} \hspace{0.5em}

This study employs modified multi-layer perceptrons (MLPs) \cite{wang2021understanding} for the branch and trunk nets in MIONets. For modeling and simulating physical systems governed by PDEs, modified MLPs have been shown to outperform conventional fully-connected neural networks (FNNs) \cite{wang2021understanding, wu2024capturing}. The key enhancement is projecting the input variables into a high-dimensional feature space by including two transformer networks. The hidden layers are then updated using a pointwise multiplication operation according to the following forward propagation rule:
\begin{equation}
    \begin{aligned}
    U & =\sigma \circ\left(XW^1 +b^1\right), V=\sigma \circ\left(XW^2+b^2\right), \\
    H^{[1]} & =\sigma \circ\left(XW^{z, 1}+b^{z, 1}\right), \\
    Z^{[k]} & =\sigma \circ\left(H^{[k]}W^{z,k} +b^{z,k}\right), \quad k=1, \ldots, K-1, \\
    H^{[l+1]} & =\left(1-Z^{[k]}\right) \odot U+Z^{[k]} \odot V, \quad k=1, \ldots, K-1, \\
    f_\theta(x) & = H^{[K]}W^{z,K}+b^{z,K}.
    \end{aligned}
\end{equation}
Here, $W^1, W^2 \in \mathbb{R}^{m_0 \times m_1}$ and $b^1, b^2 \in \mathbb{R}^{m_1}$. $W^{z,k} \in \mathbb{R}^{m_k\times m_{k+1}}$ and $b^{z,k} \in \mathbb{R}^{m_{k+1}}$ represent the weight matrix and bias vector for the $k$-th layer, respectively. The superscript $m_k$ for $k=0,\cdots,K$ denotes the dimension of the $k$-th layer, where $m_0$ and $m_K$ represent the input and output dimensions, respectively. Symbols $\sigma$, $\circ$ and $\odot$ represent scalar nonlinear activation function, element-wise operation and element-wise multiplication respectively. $X$ is a design matrix of the input data points. By incorporating residual connections and multiplicative interactions to update the hidden states, this neural network architecture can significantly improve predictive accuracy \cite{wang2021understanding}.

\subsection{Solving the two kinetic equations with PI-MIONets}

This section aims to utilize the PI-MIONets to solve the VPFP system \eqref{vlasov}-\eqref{poisson} and the semiconductor Boltzmann equation \eqref{semiB}-\eqref{poisson} with the high-field scaling. First, we establish the initial-boundary value problem corresponding to their generic model \eqref{transport}-\eqref{poisson} within a bounded domain $(t, x, v)\in\mathcal{T} \times \mathcal{D} \times \Omega$:
\begin{align}
    \partial_t f+v\cdot\nabla_x f-\frac{1}{\varepsilon}\nabla_x \phi \cdot \nabla_v f =\frac{1}{\varepsilon} \mathscr{Q}(f), &\quad(t, x, v) \in \mathcal{T} \times \mathcal{D} \times \Omega, \label{setting1} \\
    -\triangle_x \phi(t, x) =\rho(t, x)-h(x), &\quad(t, x) \in \mathcal{T} \times \mathcal{D}, \label{setting2}
\end{align}
where the operator $\mathscr{Q}$ can be either the Fokker-Planck operator $\mathscr{Q}_{\text{FP}}$ or the semiconductor Boltzmann operators $\mathscr{Q}_{\text{nond}}$ and $\mathscr{Q}_{\text{deg}}$. The initial condition (IC) on $f$ is set as 
\[
f(t=0, x, v)=f_0(x, v)
\]
where $f_0$ is a predefined function, and the initial electric potential $\phi(t=0, x)$ is determined by solving the Poisson equation \eqref{setting2}. The boundary condition (BC) on $f$ is denoted as 
\[
\mathscr{B}f(t,x,v)=F_{\mathrm{B}}, \quad {\text  {for}}\quad 
(t, x, v) \in \mathcal{T} \times \partial \mathcal{D} \times \Omega,
\]
while the BC on $\phi$ is expressed as 
\[
\mathscr{B}\phi(t,x)=\Psi_{\mathrm{B}}, \quad {\text {for}} \quad  (t, x) \in \mathcal{T} \times \partial \mathcal{D}.
\]
Here, $F_{\mathrm{B}}$ and $\Psi_{\mathrm{B}}$ are given functions, and $\partial \mathcal{D}$ represents the boundary of the domain $\mathcal{D}$. We mainly consider periodic BCs in $x$-direction in this study, the details of which will be outlined in Sect. 5.1. Under this setting, the initial distribution $f_0$ and the background charge $h$ are two variable parameters within specified ranges.

\begin{figure}[!htb]
	\centering 
	\subfigbottomskip=5pt 
	\subfigcapskip=-5pt 
    \subfigure[The architecture of the MIONet with modified MLPs for approximating the solution operator $\mathcal{F}$.]{
	\includegraphics[width=0.6\linewidth]{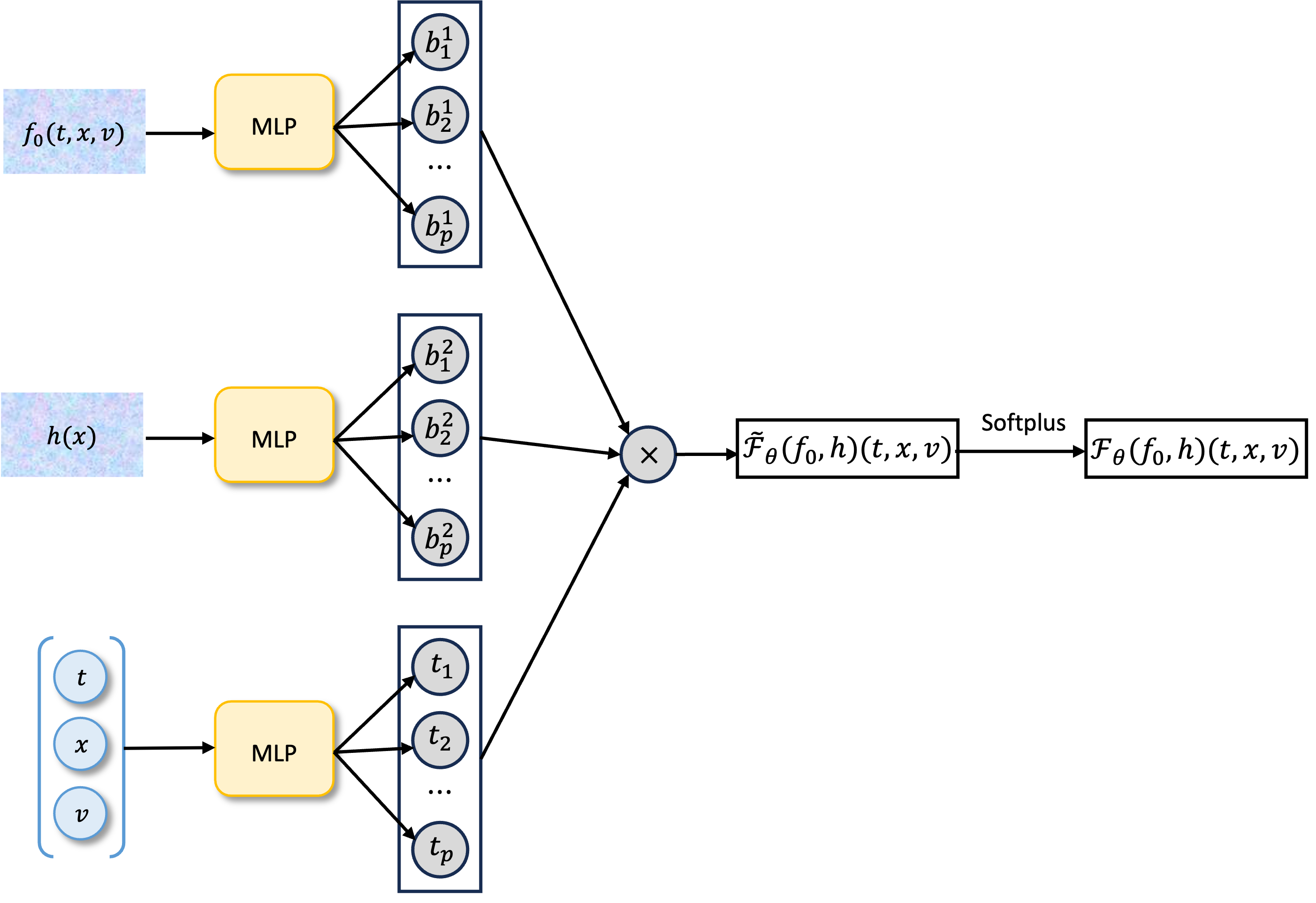}}

   \subfigure[The architecture of the PI-MIONet.]{
   \includegraphics[width=0.8\linewidth]{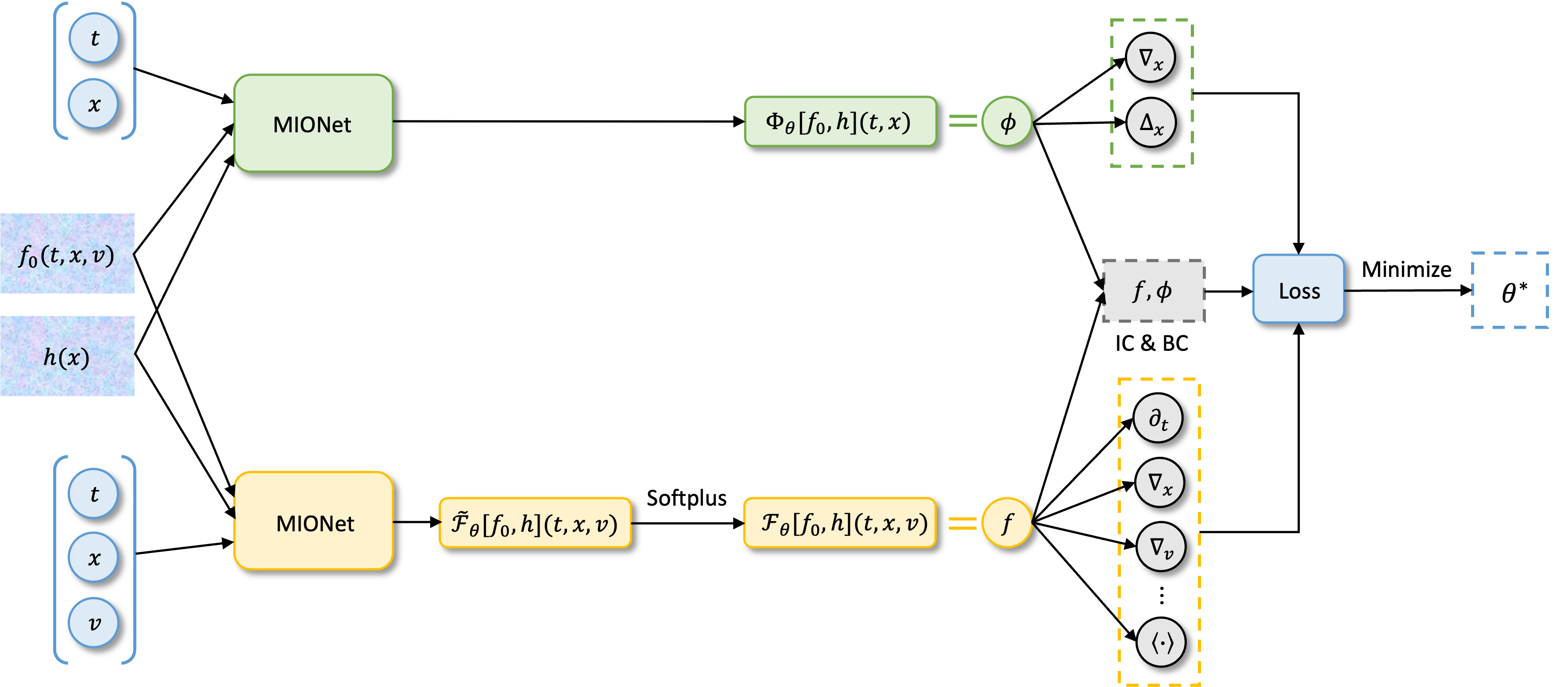}
   }
    \caption{Schematic illustrations of the physics-informed MIONet (PI-MIONet) architecture generic to both the VPFP system and the semiconductor Boltzmann equation in the high-field regime. (a) The architecture of the MIONet with modified MLPs for the parameterized operator $\mathcal{F}_\theta$ within the PI-MIONet framework. (b) The overall architecture of the PI-MIONet. }
\label{PI-MIONet2}
\end{figure}

When applying the PI-MIONet framework to the above initial-boundary value problem, our objective is to learn two solution operators that map $f_0$ and $h$ to their respective solutions, $f$ and $\phi$. These solution operators are denoted as $\tilde{\mathcal{F}}$ for $f$ and $\Phi$ for $\phi$, given by:
\begin{equation}
    \tilde{\mathcal{F}}:\left(f_0, h\right) \mapsto f; \quad \Phi:\left(f_0, h\right) \mapsto \phi.
\end{equation}
To approximate $\tilde{\mathcal{F}}$ and $\Phi$ respectively, two distinct MIONets are employed. Each MIONet consists of two branch nets and a single trunk net: the branch nets individually take the discretized $f_0$ and $h$ as inputs; and the trunk net accepts the location coordinates for evaluating the output functions. The modified MLPs serve as the network architecture for all sub-networks within these MIONets. For simplicity, a unified parameter $\theta$ is introduced to represent the collection of all trainable weights and biases across both MIONets for approximating $\tilde{\mathcal{F}}$ and $\Phi$, which are then parameterized as $\tilde{\mathcal{F}}_\theta$ and $\Phi_\theta$ respectively. To ensure the positivity of the distribution function $f$, the Softplus activation function $\sigma_{+}(\cdot):=\log (1+\exp (\cdot))$ is incorporated into the output layer of the MIONet for $\tilde{\mathcal{F}}_\theta$. A new parameterized operator $\mathcal{F}_\theta$ is introduced to represent the $\tilde{\mathcal{F}}_\theta$ composed with $\sigma_{+}$:
\begin{equation}
    \mathcal{F}_\theta[f_0, h](t, x, v):=\sigma_{+}\left(\tilde{\mathcal{F}}_\theta(f_0, h)(t, x, v)\right) \approx f(t, x, v).
\end{equation}
This network overcomes the challenge of preserving the positivity of $f$, which is a major difficulty for classical numerical methods \cite{jin2011asymptotic, lee2023oppinn}. The parameterized operator $\Phi_\theta$ is expressed as:
\begin{equation}
    \Phi_\theta[f_0, h](t, x) \approx \phi(t, x).
\end{equation}
An illustrative example of the MIONet architecture for $\mathcal{F}_\theta$ is depicted in the top panel of Fig. \ref{PI-MIONet2}.

To make the MIONets physics-informed, we encode the mean square residuals of Eqs. \eqref{setting1}-\eqref{setting2} and their  IBCs into the combined loss function of the two MIONets $\mathcal{F}_\theta$ and $\Phi_\theta$. For the sake of clarity, we introduce the notation $\boldsymbol{u}^{(i)} := (f_0^{(i)}, h^{(i)})$ for $i = 1, \cdots, M$ to represent couples of sampled input functions $f_0^{(i)}$ and $h^{(i)}$, where $M$ denotes the total number of these couples. For any couple $\boldsymbol{u}=(f_0, h)$, the loss function is formulated as:
\begin{equation}\label{PIDON}
    \begin{aligned}  \mathscr{R}_{\mathrm{PI-MIONet}}^{\varepsilon, \boldsymbol{u}} & =\frac{\mu_{1}}{|\mathcal{T} \times \mathcal{D} \times \Omega|} \int_{\mathcal{T}} \int_{\mathcal{D}} \int_{\Omega}\left|\varepsilon \partial_t \mathcal{F}_\theta+\varepsilon v \cdot\nabla_x \mathcal{F}_\theta-\nabla_x \Phi_\theta \cdot \nabla_v \mathcal{F}-\mathscr{Q}(\mathcal{F}_\theta)\right|^2 \mathrm{d}v\mathrm{d}x\mathrm{d}t \\
    & +\frac{\mu_{1}}{|\mathcal{T} \times \mathcal{D}|} \int_{\mathcal{T}} \int_{\mathcal{D}}\left|-\Delta_x \Phi_\theta-\left(\rho-h\right)\right|^2 \mathrm{d}x \mathrm{d} t \\
    & +\frac{\mu_{2}}{|\mathcal{D} \times \Omega|} \int_{\mathcal{D}} \int_{\Omega}\left|\mathcal{F}_\theta(t=0)-f_0\right|^2 \mathrm{d} v \mathrm{d} x\\
    & +\frac{\mu_{2}}{|\mathcal{D}|} \int_{\mathcal{D}}\left|-\Delta_x\Phi_\theta(t=0)-(\rho(t=0)-h)\right|^2 \mathrm{d} x\\
    & +\frac{\mu_{3}}{|\mathcal{T} \times \partial \mathcal{D} \times \Omega|} \int_{\mathcal{T}} \int_{\partial \mathcal{D}} \int_{\Omega}\left|\mathscr{B}\mathcal{F}_\theta-F_{\mathrm{B}}\right|^2 \mathrm{d} v \mathrm{d} x \mathrm{d} t \\
    & +\frac{\mu_{3}}{|\mathcal{T} \times \partial \mathcal{D}|} \int_{\mathcal{T}} \int_{\partial \mathcal{D}}\left|\mathscr{B} \Phi_\theta-\Phi_{\mathrm{B}}\right|^2 \mathrm{d} x \mathrm{d} t,
    \end{aligned}
\end{equation}
where $\mu_1$, $\mu_2$, and $\mu_3$ are weights designed to balance the contributions of different terms. To ensure a thorough exploration of function spaces, the total training loss is calculated as the average of $\mathscr{R}_{\mathrm{PI-MIONet}}^{\varepsilon, \boldsymbol{u}}$ across all samples $\boldsymbol{u}\in\left\{\boldsymbol{u}^{(i)}\right\}_{i=1}^M$, and defined as:
\begin{equation}
    \mathscr{R}_{\text {PI-MIONet}}^{\varepsilon}=\frac{1}{M} \sum_{i=1}^M \mathscr{R}_{\text {PI-MIONet}}^{\varepsilon, \boldsymbol{u}^{(i)}}.
\end{equation}
The bottom panel in Fig. \ref{PI-MIONet2} depicts a schematic diagram of this PI-MIONet architecture generic to both the VPFP system and the semiconductor Boltzmann equation.

\subsection{The limitation of PI-MIONets with a small scale parameter}\label{3.3}

This section presents the limitation of PI-MIONets in accurately capturing the equilibrium solutions when the scale parameter is small. Specifically, we check whether the loss function of the PI-MIONet can (formally) converge to the loss of the correct high-field limit equations as $\varepsilon\rightarrow0$. This desired property of the loss function exactly corresponds to the asymptotic-preserving (AP) mechanism. Our analysis focuses on the residual terms associated with Eqs. \eqref{setting1}-\eqref{setting2} in the loss function $\mathscr{R}_{\mathrm{PI-MIONet}}^{\varepsilon, \boldsymbol{u}}$: 
\begin{equation}\label{PDEs1}
    \begin{aligned}
    \mathscr{R}_{\text {PDEs}}^{\varepsilon, \boldsymbol{u}} & =\frac{\mu_{1}}{|\mathcal{T} \times \mathcal{D} \times \Omega|} \int_{\mathcal{T}} \int_{\mathcal{D}} \int_{\Omega}\left|\varepsilon \partial_t \mathcal{F}_\theta+\varepsilon v \cdot \nabla_x \mathcal{F}_\theta-\nabla_x \Phi_\theta \cdot \nabla_v \mathcal{F}-\mathscr{Q}\left(\mathcal{F}_\theta\right)\right|^2 \mathrm{d} v \mathrm{d} x \mathrm{d} t \\
    & \,\,\,\,\,\,  +\frac{\mu_{1}}{|\mathcal{T} \times \mathcal{D}|} \int_{\mathcal{T}} \int_{\mathcal{D}}\left|-\Delta_x \Phi_\theta-(\rho-h)\right|^2 \mathrm{d} x \mathrm{d} t.
    \end{aligned}
\end{equation}
Taking $\varepsilon \rightarrow 0$ in this loss, it formally leads to the limiting loss:
\begin{equation}\label{PDEs2}
    \begin{aligned}
    \mathscr{R}_{\text {PDEs}}^{\boldsymbol{u}} & =\frac{\mu_{1}}{|\mathcal{T} \times \mathcal{D} \times \Omega|} \int_{\mathcal{T}} \int_{\mathcal{D}} \int_{\Omega}\left|\nabla_x \Phi_\theta \cdot \nabla_v \mathcal{F}+\mathscr{Q}\left(\mathcal{F}_\theta\right)\right|^2 \mathrm{d} v \mathrm{d} x \mathrm{d} t \\
    &  \,\,\,\,\,\, +\frac{\mu_{1}}{|\mathcal{T} \times \mathcal{D}|} \int_{\mathcal{T}} \int_{\mathcal{D}}\left|-\Delta_x \Phi_\theta-(\rho-h)\right|^2 \mathrm{d} x \mathrm{d} t.
    \end{aligned}
\end{equation}
This limiting loss is essentially the mean squared residuals of the subsequent system:
\begin{equation}\label{PDEs3}
    \left\{\begin{aligned}
    -\nabla_x\phi \cdot \nabla_v f & = \mathscr{Q}(f), \\
    -\Delta_x \phi & =\rho-h. \\
    \end{aligned}\right.
\end{equation}
However, this limiting system does not align with the expected high-field limit equations as presented in Eq. \eqref{newvfp}, Eq. \eqref{nond_eq1} or Eq. \eqref{deg_limit}. One possible reason is that at small values of $\varepsilon$, the loss $\mathscr{R}_{\text {PI-MIONet}}^{\varepsilon}$ is primarily dominated by its leading-order terms, thereby only capturing the single scale behavior. Therefore, training the PI-MIONet using the loss function $\mathscr{R}_{\text {PI-MIONet}}^{\varepsilon}$ may suffer from significant inaccuracy  issue in capturing the desired macroscopic limit behavior. This instability will be further verified through numerical experiments in Sect. 5.

\section{Methodology}

In this section, we propose a new asymptotic-preserving multiple-input DeepONet (AP-MIONet) to address the multiscale challenge discussed in Sect. \ref{3.3}. The pivotal procedure of the AP-MIONet method is to design a loss function with the asymptotic-preserving (AP) property. A major challenge in designing an AP loss for the high-field regime arises from the frequent absence of an explicit expression for the local equilibrium, a main component of the AP loss design. Our method begins by constructing a new model that incorporates the mass conservation law with the original kinetic system, generic to both the VPFP and semiconductor Boltzmann equations. We then reformulate the loss function of the PI-MIONet by encoding the residuals of the proposed model, thereby establishing the AP-MIONet. We conclude by formally validating the AP property of the  established loss function, demonstrating the uniform stability of the AP-MIONet with respect to the small scale parameter.

\subsection{Enforcing the mass conservation law}

In the high-field regime, the explicit local equilibrium is often absent, particularly for the semiconductor Boltzmann equation. Hence we cannot leverage the micro-macro decomposition \cite{jin2023asymptotic1, lu2022solving, li2022model}, which is efficient for AP loss design in the diffusion regime, nor can we utilize the even-odd decomposition \cite{jin2023asymptotic, wu2024capturing} because a ``non-stiff'' force term cannot be derived. Here, we adopt the idea of enforcing the mass conservation law besides the original kinetic system to develop a new model for AP loss design, inspired by the classical AP schemes for the high-field regime \cite{jin2011asymptotic, jin2013asymptotic}. Let us begin with the original kinetic system \eqref{transport}-\eqref{poisson}:
\begin{subequations}
    \begin{equation}\tag{1a}
        \partial_t f+v\cdot\nabla_x f-\frac{1}{\varepsilon}\nabla_x \phi \cdot \nabla_v f =\frac{1}{\varepsilon}\mathscr{Q}(f),
    \end{equation}
    \begin{equation}\tag{1b}
        -\triangle_x \phi(t, x) =\rho(t, x)-h(x),
    \end{equation}
\end{subequations} 
where the operator $\mathscr{Q}$ can be either the Fokker-Planck operator $\mathscr{Q}_{\text{FP}}$ or the semiconductor Boltzmann operator $\mathscr{Q}_{\text{semiB}}$. One integrates the first equation \eqref{transport} over $v$ and obtains the mass conservation law:
\begin{equation}\label{law3}
    \partial_t \rho+\nabla_x \cdot\langle v f\rangle=0.
\end{equation}
The new model incorporates the original system Eq. \eqref{transport}-Eq. \eqref{poisson}, the mass conservation law \eqref{law3} and the conservation condition $\rho=\langle f\rangle$, established as follows:
\begin{equation}\label{law4}
    \left\{\begin{aligned}
        &\partial_t f+v\cdot\nabla_x f-\frac{1}{\varepsilon}\nabla_x \phi \cdot \nabla_v f =\frac{1}{\varepsilon}\mathscr{Q}(f), \\
        &\partial_t \rho+\nabla_x \cdot\langle v f\rangle  =0, \\
        &-\triangle_x \phi  =\rho-h, \\
        &\rho =\langle f\rangle .
    \end{aligned}\right.
\end{equation}
This model is equivalent to the original system, since the redundant equation \eqref{law3} is merely the integration of Eq.\eqref{transport} over the velocity space, but it is an easy yet crucial way to guarantee the conservation of mass when designing the loss function. 

We now check whether the correct high-field limit equations can be derived from the new model Eqs. \eqref{law4} with the Fokker-Planck and semiconductor Boltzmann collision operators. These derivations closely follow the procedures outlined in Sect. 3.1 for the VPFP system and in Sect. 3.2 for the semiconductor Boltzmann equation, respectively. First, one considers the Fokker-Planck operator $\mathscr{Q}_\text{FP}$. One multiplies the first equation in Eqs. \eqref{law4} by velocity $v$ and integrates it over the velocity space $\mathbb{R}^N$. On taking the limit as $\varepsilon \rightarrow 0$, one can infer $j = -\rho(\nabla_x \phi)$. Substituting this inference into the second equation in Eqs. \eqref{law4}, the limit system degenerates to: 
\begin{equation}
    \left\{\begin{array}{l}
    \partial_t \rho-\nabla_x \cdot\left(\rho \nabla_x \phi\right)=0, \\
    -\Delta_x \phi=\rho-h(x),
    \end{array}\right.
\end{equation}
which aligns with the limit equation \eqref{limit} presented in Sect. 3.1. 

We next focus on the semiconductor Boltzmann operators $\mathscr{Q}_\text{nond}$ and $\mathscr{Q}_\text{deg}$. For the non-degenerate case, as $\varepsilon\rightarrow0$ in the first equation in Eqs.\eqref{law4}, the distribution $f$ converges to $\rho(t, x) F_{E(t, x)}(v)$, where $F_{E(t, x)}(v)$ complies with Eq.\eqref{nond_eq2}. Substituting $f=\rho(t, x) F_{E(t, x)}(v)$ into the second equation of Eqs. \eqref{law4}, one derives the associated limit equation Eq.\eqref{deriv5}. For the degenerate case, when sending $\varepsilon$ to zero in the first equation in Eqs. \eqref{law4}, $f$ converges to $F(\rho(t, x), E(t, x))(v)$, where $F(\rho(t, x), E(t, x))(v)$ satisfies Eq. \eqref{deriv7}. Incorporating $f=F(\rho(t, x), E(t, x))(v)$ into the second equation of Eqs. \eqref{law4}, one deduces the limit equation Eq. \eqref{deg_limit}.

The proposed model \eqref{law4} presents several advantages that can be leveraged to construct an AP loss. Foremost, it can accurately derive the high-field limit equation without requiring an explicit local equilibrium. Therefore, it provides a generic strategy to develop an AP-MIONet method compatible with both the VPFP and semiconductor Boltzmann equations, details of which will be elaborated in the following section. Secondly, this model preserves the mass and handles non-equilibrium states, tackling two difficulties commonly faced by neural networks \cite{jin2023asymptotic2, jin2023asymptotic}. Specifically, the mass conservation equation naturally guarantees the preservation of mass,  while the kinetic equation for the probability density distribution captures the dynamics of far-from-equilibrium states.

\subsection{Asymptotic-preserving multiple-input DeepONets}

We present the definition of the asymptotic-preserving multiple-input DeepONets (AP-MIONets) as follows:

\begin{definition}{Asymptotic-preserving multiple-input DeepONets (AP-MIONets)}\\
For a multiscale equation with multiple variable parameters, assume that its solution operator is parameterized by a multiple-input DeepONet (MIONet). This MIONet is trained by minimizing a loss function that incorporates the residual of the governing equation via a gradient-based optimization method. A MIONet is called an Asymptotic-preserving multiple-input DeepONet (AP-MIONet) if, as the physical scale parameter approaches zero, the loss function of the microscopic equation automatically transitions to the loss of the corresponding macroscopic equation.
\end{definition}

 In other words, the key idea behind developing an AP-MIONet is to design a loss function with the AP property, thereby ensuring that the MIONet can preserve the correct equilibrium solution as the scale parameter decreases.

\begin{figure}[!htb]
	\centering  
	\includegraphics[width=0.8\linewidth]{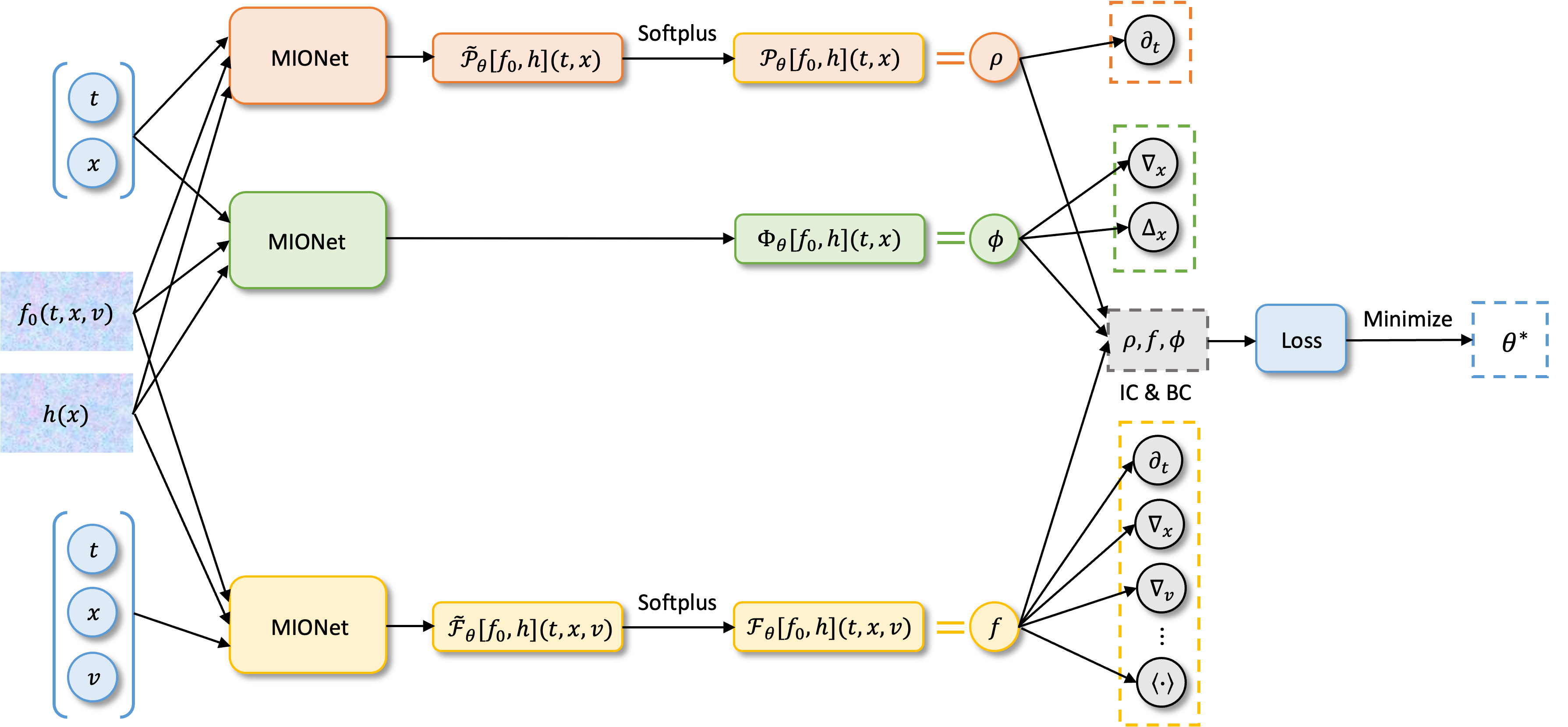}
	\caption{The architecture of the Asymptotic-Preserving Multiple-Input DeepONet (AP-MIONet) generic to both the VPFP system and the semiconductor Boltzmann equation in the high-field regime.}
\label{AP-MIONet}
\end{figure}

We now provide the details of the proposed AP-MIONet method. To design a loss function with the AP property, the PI-MIONet is employed to solve the new model Eqs.\eqref{law4}, from which the high-field limit can be readily derived, instead of the original system Eq. \eqref{transport}-Eq. \eqref{poisson}. First, we define the initial-boundary conditions for this model, which are specified within the bounded domain $(t,x,v)\in \mathcal{T} \times \mathcal{D} \times \Omega$. The ICs on $f$ is $f(t=0, x, v)=f_0(x, v)$ where $f_0$ is a given function. The initial mass density $\rho(t=0,x)$ is calculated using the relation $\rho(t=0,x)=\langle f(t=0, x, v)\rangle$, and the initial electric potential $\phi(t=0, x)$ is determined by the Poisson equation \eqref{setting2}. The BCs on $f, \rho$ and $\phi$ are denoted as $\mathscr{B} f=F_B, \mathscr{B}\rho=P_B$ and $\mathscr{B} \phi=\Psi_B$ respectively, where $F_B, P_B$ and $\Psi_B$ are given functions. We primarily consider periodic BCs in the $x$-direction for $f,\rho$ and $\phi$, with further details provided in Sect. 5.1. When solving this initial-boundary value problem by the PI-MIONet method, our objective is to learn three solution operators that map from $f_0$ and $h$ to their respective solutions $f$, $\rho$, and $\phi$, denoted as $\tilde{\mathcal{F}}$, $\tilde{\mathcal{P}}$, and $\Phi$:
\begin{equation}
    \tilde{\mathcal{F}}:\left(f_0, h\right) \mapsto f, \quad\tilde{\mathcal{P}}:\left(f_0, h\right) \mapsto \rho, \quad\Phi:\left(f_0, h\right) \mapsto \phi.
\end{equation}
Three distinct MIONets with modified MLP are employed to approximate the solution operators $\tilde{\mathcal{F}}, \tilde{\mathcal{P}}$ and $\Phi$. Each MIONet comprises two branch nets and a single trunk net. To preserve the positivity of the distribution $f$ and density $\rho$, two operators $\mathcal{F}_\eta$ and $\mathcal{P}_\eta$ are introduced to represent the parameterized operators $\tilde{\mathcal{F}}_\eta$ and $\tilde{\mathcal{P}}_\eta$ with the Softplus activation function $\sigma_{+}$ applied. Here, $\eta$ denotes
all trainable network parameters for brevity. These operators are defined as:
\begin{subequations}
    \begin{equation}
        \mathcal{F}_\eta[f_0, h](t, x, v):=\sigma_{+}\left(\tilde{\mathcal{F}}_\eta(f_0, h)(t, x, v)\right) \approx f(t, x, v),
    \end{equation}
    \begin{equation}
        \mathcal{P}_\eta[f_0, h](t, x):=\sigma_{+}\left(\tilde{\mathcal{P}}_\eta(f_0, h)(t, x)\right) \approx\rho(t, x).
    \end{equation}
\end{subequations}
The parameterized operator $\Phi_\eta[f_0, h](t, x)$ is formulated as:
\begin{equation}
    \Phi_\eta[f_0, h](t, x)\approx \phi(t, x).
\end{equation}
The loss function of the AP-MIONet, for each input couple $\boldsymbol{u}=(f_0, h)$, is formulated by taking the mean squared residuals of the model Eqs. \eqref{law4} and its associated IBCs as follows:
\begin{equation}\label{APDON}
    \begin{aligned}
    \mathscr{R}_{\mathrm{AP-MIONet}}^{\varepsilon, \boldsymbol{u}} & =\frac{\lambda_{1}}{|\mathcal{T} \times \mathcal{D}|} \int_{\mathcal{T}} \int_{\mathcal{D}}\left|\partial_t \mathcal{P}_\eta+\nabla_x\cdot\left\langle v \mathcal{F}_\eta\right\rangle\right|^2 \mathrm{d} x \mathrm{d} t\\
    &\frac{\lambda_{1}}{|\mathcal{T} \times \mathcal{D} \times \Omega|} \int_{\mathcal{T}} \int_{\mathcal{D}} \int_{\Omega}\left|\varepsilon \partial_t \mathcal{F}_\eta+\varepsilon v \cdot\nabla_x \mathcal{F}_\eta-\nabla_x \Phi_\eta \cdot \nabla_v \mathcal{F}-\mathscr{Q}(\mathcal{F}_\eta)\right|^2 \mathrm{d} v \mathrm{d} x \mathrm{d}t \\
    & +\frac{\lambda_{1}}{|\mathcal{T} \times \mathcal{D}|} \int_{\mathcal{T}} \int_{\mathcal{D}}\left|-\Delta_x \Phi_\eta-\left(\mathcal{P}_\eta-h\right)\right|^2 \mathrm{d} x\mathrm{d} t\\ &+\frac{\lambda_{2}}{|\mathcal{T} \times \mathcal{D}|} \int_{\mathcal{T}} \int_{\mathcal{D}}\left|\left\langle \mathcal{F}_\eta\right\rangle-\mathcal{P}_\eta\right|^2 \mathrm{d} x\mathrm{d} t\\
    & +\frac{\lambda_{3}}{|\mathcal{D}|} \int_{\mathcal{D}} \left|\mathcal{P}_\eta(t=0)-\langle\mathcal{F}_\eta(t=0)\rangle\right|^2 \mathrm{d} x +\frac{\lambda_{3}}{|\mathcal{D} \times \Omega|} \int_{\mathcal{D}} \int_{\Omega}\left|\mathcal{F}_\eta(t=0)-f_0\right|^2 \mathrm{d} v \mathrm{d} x\\
    & +\frac{\lambda_{3}}{|\mathcal{D}|} \int_{\mathcal{D}}\left|-\Delta_x\Phi_\eta(t=0)-(\mathcal{P}_\eta(t=0)-h)\right|^2 \mathrm{d} x\\
    & +\frac{\lambda_{4}}{|\mathcal{T} \times \partial \mathcal{D}|} \int_{\mathcal{T}} \int_{\partial \mathcal{D}}\left|\mathscr{B} \mathcal{P}_\eta-P_{\mathrm{B}}\right|^2 \mathrm{d} x \mathrm{d} t  +\frac{\lambda_{4}}{|\mathcal{T} \times \partial \mathcal{D} \times \Omega|} \int_{\mathcal{T}} \int_{\partial \mathcal{D}} \int_{\Omega}\left|\mathscr{B}\mathcal{F}_\eta-F_{\mathrm{B}}\right|^2 \mathrm{d} v \mathrm{d} x \mathrm{d} t \\
    & +\frac{\lambda_{4}}{|\mathcal{T} \times \partial \mathcal{D}|} \int_{\mathcal{T}} \int_{\partial \mathcal{D}}\left|\mathscr{B} \Phi_\eta-\Phi_{\mathrm{B}}\right|^2 \mathrm{d} x \mathrm{d} t, \\
    \end{aligned}
\end{equation}
where $\lambda_1, \lambda_2, \lambda_3$ and $\lambda_4$ are penalty weights that can be fine-tuned. The total loss is defined as the average of $\mathscr{L}_{\text {AP-MIONet}}^{\varepsilon, \boldsymbol{u}^{(i)}}$ across all input samples $\boldsymbol{u} \in\left\{\boldsymbol{u}^{(i)}\right\}_{i=1}^M$:
\begin{equation}\label{APDON1}
    \mathscr{R}_{\text {AP-MIONet }}^{\varepsilon}=\frac{1}{M} \sum_{i=1}^M \mathscr{R}_{\text {AP-MIONet}}^{\varepsilon, \boldsymbol{u}^{(i)}}.
\end{equation}
Fig. \ref{AP-MIONet} depicts a schematic illustration of the AP-MIONet architecture.

We then formally validate the AP property of the loss function \eqref{APDON1}. It suffices to examine the residual terms of the model Eqs.\eqref{law4} within the loss \eqref{APDON}:
\begin{equation}\label{APDON2}
    \begin{aligned}    \mathscr{R}_{\mathrm{Model}}^{\varepsilon, \boldsymbol{u}} & =\frac{\lambda_{1}}{|\mathcal{T} \times \mathcal{D}|} \int_{\mathcal{T}} \int_{\mathcal{D}}\left|\partial_t \mathcal{P}_\eta+\nabla_x\cdot\left\langle v \mathcal{F}_\eta\right\rangle\right|^2 \mathrm{d} x \mathrm{d} t\\
    &\frac{\lambda_1}{|\mathcal{T} \times \mathcal{D} \times \Omega|} \int_{\mathcal{T}} \int_{\mathcal{D}} \int_{\Omega}\left|\varepsilon \partial_t \mathcal{F}_\eta+\varepsilon v \cdot\nabla_x \mathcal{F}_\eta-\nabla_x \Phi_\eta \cdot \nabla_v \mathcal{F}-\mathscr{Q}(\mathcal{F}_\eta)\right|^2 \mathrm{d} v \mathrm{d} x \mathrm{d}t \\
    & +\frac{\lambda_1}{|\mathcal{T} \times \mathcal{D}|} \int_{\mathcal{T}} \int_{\mathcal{D}}\left|-\Delta_x \Phi_\eta-\left(\mathcal{P}_\eta-h\right)\right|^2 \mathrm{d} x\mathrm{d} t \\
    &+\frac{\lambda_2}{|\mathcal{T} \times \mathcal{D}|} \int_{\mathcal{T}} \int_{\mathcal{D}}\left|\left\langle \mathcal{F}_\eta\right\rangle-\mathcal{P}_\eta\right|^2 \mathrm{d} x\mathrm{d} t.\\
    \end{aligned}
\end{equation}
On taking the limit $\varepsilon \rightarrow 0$ in this loss, one formally obtains the limiting loss:
\begin{equation}
    \begin{aligned}    \mathscr{R}_{\mathrm{Model}}^{\boldsymbol{u}} & =\frac{\lambda_1}{|\mathcal{T} \times \mathcal{D}|} \int_{\mathcal{T}} \int_{\mathcal{D}}\left|\partial_t \mathcal{P}_\eta+\nabla_x\cdot\left\langle v \mathcal{F}_\eta\right\rangle\right|^2 \mathrm{d} x \mathrm{d} t\\
    &\frac{\lambda_1}{|\mathcal{T} \times \mathcal{D} \times \Omega|} \int_{\mathcal{T}} \int_{\mathcal{D}} \int_{\Omega}\left|\nabla_x \Phi_\eta \cdot \nabla_v \mathcal{F}+\mathscr{Q}(\mathcal{F}_\eta)\right|^2 \mathrm{d} v \mathrm{d} x \mathrm{d}t \\
    & +\frac{\lambda_1}{|\mathcal{T} \times \mathcal{D}|} \int_{\mathcal{T}} \int_{\mathcal{D}}\left|-\Delta_x \Phi_\eta-\left(\mathcal{P}_\eta-h\right)\right|^2 \mathrm{d} x\mathrm{d} t \\
    &+\frac{\lambda_2}{|\mathcal{T} \times \mathcal{D}|} \int_{\mathcal{T}} \int_{\mathcal{D}}\left|\left\langle \mathcal{F}_\eta\right\rangle-\mathcal{P}_\eta\right|^2 \mathrm{d} x\mathrm{d} t.\\
    \end{aligned}
\end{equation}
This limiting loss is essentially the mean squared residuals of the following system:
\begin{equation}
    \left\{\begin{aligned}
    \partial_t\rho+\nabla_x\cdot\langle vf\rangle &= 0 ,\\
    -\nabla_x\phi \cdot \nabla_v f & = \mathscr{Q}(f), \\
    -\Delta_x \phi & =\rho-h, \\
    \rho &= \langle f \rangle.
    \end{aligned}\right.
\end{equation}
By revisiting and replicating the derivation procedures detailed in Sect. 4.1 for the above system, we can recover the correct high-field limit equations for the VPFP system and the semicondutor Boltzmann equation. Therefore, the AP property of the proposed AP-MIONet loss function is  verified.

Our approach shows novelty in successfully deriving the AP property (at least formally) for the semiconductor Boltzmann equation in its degenerate form. In this case, the nonlinearity of the collision operator poses a significant challenge in theoretically verifying the AP property of numerical schemes \cite{jin2013asymptotic}. We tackle this challenge by leveraging theoretical proofs from the governing continuous equation, thus establishing the AP property for the proposed loss function.

\section{Numerical results}

In this section, we present comprehensive numerical results for a diverse array of test problems across the kinetic and high-field regimes. These results demonstrate the efficacy and versatility of the proposed AP-MIONet in accurately capturing the complex multiscale dynamics.

\subsection{Experiment setup}

This part provides detailed settings for the practical implementation of the numerical experiments. To generate the dataset of input functions, we sample $M=640$ discretized input function couples, denoted as $\{\boldsymbol{u}^{(i)}\}_{i=1}^{M}=\{(f_0^{(i)}, h^{(i)})\}_{i=1}^{M}$. Each $f_0^{(i)}$ represents the initial condition $f_0(x,v)$ evaluated at a collection of fixed sensors $\{(x_l,v_q)\}_{1 \leq l \leq L, 1 \leq q \leq Q}$, where $\{x_l\}_{l=1}^L, \{v_q\}_{q=1}^{Q}$ are equidistant grid points in the $x$ and $v$ directions respectively. Each background charge $h^{(i)}$ denotes $ \left[h^{(i)}(x_1), h^{(i)}(x_2), \ldots, h^{(i)}(x_L)\right]$, assessed at equi-spaced locations $\{x_l\}_{l=1}^L$ along the $x$-axis. Here, both $L$ and $Q$ are set to $32$. These samples are divided into training and testing datasets at a $4:1$ ratio, leading to $M_\text{train}=512$ training couples and $M_\text{test}=128$ testing couples. 

In the MIONet architecture, each branch and trunk net employs a modified MLP consisting of $5$ hidden layers ($K = 5$), with each layer containing $64$ neurons. The $swish$ activation function \cite{ramachandran2017searching} is utilized in all hidden layers due to its stability in training dynamics and superior performance compared to classical activation functions. These networks are initialized using the Glorot normal scheme \cite{glorot2010understanding}, which effectively stabilizes the training process and accelerates convergence. Periodic BCs in the $x$-direction are applied to the solutions $\rho, f$ and $\phi$ in all numerical examples. To improve the numerical performance, the periodic BCs are directly encoded into the MIONets as hard constraints by modifying network architectures \cite{lu2022comprehensive}. Specifically, for the trunk net in each MIONet, the spatial input $x$ is transformed using  Fourier basis functions: 
\begin{equation}
    \{\cos (\omega x), \sin (\omega x), \cos (2 \omega x), \sin (2 \omega x), \cdots\},
\end{equation}
where $\omega = 2 \pi/P$ and $P$ denotes the period. The inherent periodicity of these basis functions ensures that the output of the MIONets exhibits the desired periodic behavior in $x$ \cite{dong2021method}. In this study, the Fourier feature set is restricted to the first two terms -- ${\cos (\omega x), \sin (\omega x)}$ -- to balance model complexity and computational efficiency.

In training AP-MIONets, the integral-based loss functions are approximated using the Monte Carlo method, leading to the \textit{empirical loss functions}. To ensure effective training, we embed the exact solutions of the initial conditions (ICs) for $f$, $\rho$, and $\phi$ into the loss function. Specifically, the original ICs for $\rho$, $f$ and $\phi$:
\begin{equation}\label{ICs}
f(t=0) = f_0(x, v),\quad\rho(t=0)=\langle f(t=0)\rangle,\quad -\Delta_x\phi(t=0) = \rho(t=0)-h,
\end{equation}
can be reformulated as:
\begin{equation}
    f(t=0) = f_0(x, v),\quad\rho(t=0)=\rho_0(x),\quad \phi(t=0) = \phi_0(x),
\end{equation}
where $\rho_0$ and $\phi_0$ are derived from solving Eqs.\eqref{ICs}.
For cases lacking analytical expressions for $\rho(t=0)$, $f(t=0)$ and $\phi(t=0)$, we can numerically generate data for these conditions and incorporate them as regularization terms into the loss function. The empirical loss function for the AP-MIONet is then constructed as: 
\begin{equation}\label{empirical}
    \begin{aligned}        \mathscr{L}^\varepsilon_{\text{AP-MIONet}} &:= \frac{\lambda_1}{M_{\text{train}}N_{\text{dom}}}\sum_{i=1}^{M_{\text{train}}}\sum_{j=1}^{N_{\text{dom}}}\left| ( \varepsilon\partial_t\mathcal{F}_\eta -\varepsilon v\cdot\nabla_x\mathcal{F}_\eta - \nabla_x\Phi_\eta \cdot\nabla_v\mathcal{F}_\eta -\mathscr{Q}(\mathcal{F}_\eta))(\boldsymbol{u}^{(i)})(t_{r,j}^{(i)}, x_{r,j}^{(i)}, v_{r,j}^{(i)})   \right|^2\\
    &+ \frac{\lambda_1}{M_{\text{train}}N_{\text{dom}}}\sum_{i=1}^{M_{\text{train}}}\sum_{j=1}^{N_{\text{dom}}}\left|(\partial_t\mathcal{P}_\eta - \nabla_x\cdot\langle v\mathcal{F}_\eta \rangle)(\boldsymbol{u}^{(i)})(t_{r,j}^{(i)}, x_{r,j}^{(i)})  \right|^2\\
    &+ \frac{\lambda_1}{M_{\text{train}}N_{\text{dom}}}\sum_{i=1}^{M_{\text{train}}}\sum_{j=1}^{N_{\text{dom}}}\left| (\nabla_x\Phi_\eta-(\mathcal{P}_\eta-h))(\boldsymbol{u}^{(i)})(t_{r,j}^{(i)}, x_{r,j}^{(i)})  \right|^2\\
    &+ \frac{\lambda_2}{M_{\text{train}}N_{\text{dom}}}\sum_{i=1}^{M_{\text{train}}}\sum_{j=1}^{N_{\text{dom}}}\left|(\mathcal{P}_\eta-\langle \mathcal{F}_\eta\rangle)(\boldsymbol{u}^{(i)})(t_{r,j}^{(i)}, x_{r,j}^{(i)}) \right|^2 \\
    &+ \frac{\lambda_{3,\rho}}{M_{\text{train}}N_{\text{init}}}\sum_{i=1}^{M_{\text{train}}}\sum_{j=1}^{N_{\text{init}}}\left| \mathcal{P}_\eta(\boldsymbol{u}^{(i)})(0, x_{ic,j}^{(i)})-\rho_0(x_{ic,j}^{(i)})  \right|^2 \\
    &+ \frac{\lambda_{3,f}}{M_{\text{train}}N_{\text{init}}}\sum_{i=1}^{M_{\text{train}}}\sum_{j=1}^{N_{\text{init}}}\left| \mathcal{F}_\eta(\boldsymbol{u}^{(i)})(0, x_{ic,j}^{(i)}, v_{ic,j}^{(i)})-f_0(x_{ic,j}^{(i)}, v_{ic,j}^{(i)})   \right|^2 \\
    &+ \frac{\lambda_{3,\phi}}{M_{\text{train}}N_{\text{init}}}\sum_{i=1}^{M_{\text{train}}}\sum_{j=1}^{N_{\text{init}}}\left| \Phi_\eta(\boldsymbol{u}^{(i)})(0, x_{ic,j}^{(i)})-\phi_0(x_{ic,j}^{(i)})   \right|^2. \\
    \end{aligned}
\end{equation}
Here, for each $\boldsymbol{u}^{(i)}$, $\{(t_{r,j}^{(i)}, x_{r,j}^{(i)}) \}_{j=1}^{N_\text{dom}}$ and $\{(t_{r,j}^{(i)}, x_{r,j}^{(i)}, v_{r,j}^{(i)}) \}_{j=1}^{N_{\text{dom}}}$ are uniformly sampled within the computational domain to enforce the PDE residuals, while $\{( x_{ic,j}^{(i)})\}_{j=1}^{N_{\text{init}}}$ and $\{( x_{ic,j}^{(i)}, v_{ic,j}^{(i)})\}_{j=1}^{N_{\text{init}}}$ are two sets of collocation points to impose the ICs. The number of sample points is set to $N_{\text{dom}}=2^{10}$ within the domain and $N_{\text{init}}=2^{9}$ for the ICs. In the training dataset, data points are organized into structured tuples, such as $(\boldsymbol{u}^{(i)}, (t_{r,j}^{(i)}, x_{r,j}^{(i)}))$, which represents a typical data point within the domain. For a detailed explanation of this data organization, please refer to \cite{wang2021learning}. The AP-MIONets are trained by minimizing the empirical loss function Eq. \eqref{empirical} using the Adam optimizer \cite{kingma2014adam}. The learning rate of the optimizer is initialized at $10^{-3}$ and decays with a rate of $0.9$ every $1,000$ training iterations. The batch sizes for the sample points are set to $4\times 10^4$ and $2\times 10^4$ for the domain and ICs respectively. The integration terms involving the operator $\langle\cdot\rangle$ are approximated by the Gauss-Legendre quadrature rule with $16$ nodes. All the hyper-parameters have been tuned to achieve optimal performance. The detailed convergence patterns of the training losses for all examples are provided in Appendix A.

To evaluate the performance of the AP-MIONet method, we compute the $\ell^2$-relative error of macroscopic quantities, namely the density $\rho$ and the electric field $E(t, x)$, between the reference solutions and the solutions obtained by the AP-MIONet method. Specifically, the $\ell^2$-relative error for $\rho$, calculated at uniform grid points $\left\{\left(t_i, x_j\right)\right\}_{1 \leq i \leq I, 1 \leq j \leq H}$, is defined as:
\begin{equation}
\ell^2(\rho):=\sqrt{\frac{\sum_{i, j}\left|\rho_\theta\left(t_i, x_j\right)-\rho_{\mathrm{ref}}\left(t_i, x_j\right)\right|^2}{\sum_{i, j}\left|\rho_{\mathrm{ref}}\left(t_i, x_j\right)\right|^2}},
\end{equation}
where $\rho_{\text{ref}}$ denotes the reference solution and $\rho_\theta$ represents the solution obtained by the AP-MIONet. Here, the reference solutions for the VPFP system are derived using the finite difference method detailed in \cite{jin2011asymptotic}, and those for the semiconductor Boltzmann equations are obtained using the numerical method presented in \cite{jin2013asymptotic}.

\subsection{Results for the VPFP system}

\noindent $\bullet$ \textbf{Landau damping} \hspace{0.5em} 

Landau damping is a common phenomenon in plasma, often used as a benchmark to assess numerical methods \cite{mouhot2011landau, stix1962theory}. We employ it to validate the AP property for the loss function of the AP-MIONet. The initial distribution function takes the form:
\begin{equation}
    f_0(x, v)=\frac{h}{\sqrt{2 \pi}}\exp \left(-\frac{v^2}{2}\right)(1+\alpha\cos(kx)), \quad(x, v) \in[0,2 \pi / k] \times \mathbb{R},
\end{equation}
where the wave number $k$ is set to $0.5$, and $\alpha$ is a perturbation parameter sampled uniformly from the interval $[0.04, 0.06]$. We consider a uniform background charge $h$, which remains constant in space and is sampled from $[0.9, 1.1]$. The computational domain is $[0,2 \pi / k] \times\left[v_{\min }, v_{\max }\right]$ with $-v_{\min }=v_{\max }=6$. Periodic BCs are applied in the $x$-direction for $f$ and $\phi$. This periodic BC for $f$ is directly translated into the settings for $\rho$. The ICs for the AP-MIONet method are specified as follows: the initial density $\rho_0$ is analytically derived from $\rho_0 = \langle f_0 \rangle$ yielding $\rho_0(x) = h(1 + \alpha \cos(kx))$, and the initial electric field $\phi_0$ is the solution to the Poisson equation, given by $\phi_0(x) = h\alpha\cos(kx)/k^2$.

\begin{figure}[htbp]
	\centering 
	\subfigbottomskip=5pt 
	\subfigcapskip=-5pt 
	\subfigure[Density $\rho$ (top) and electric field $E$ (bottom) in the kinetic regime $(\varepsilon=1)$.]{
	\includegraphics[width=0.7\linewidth]{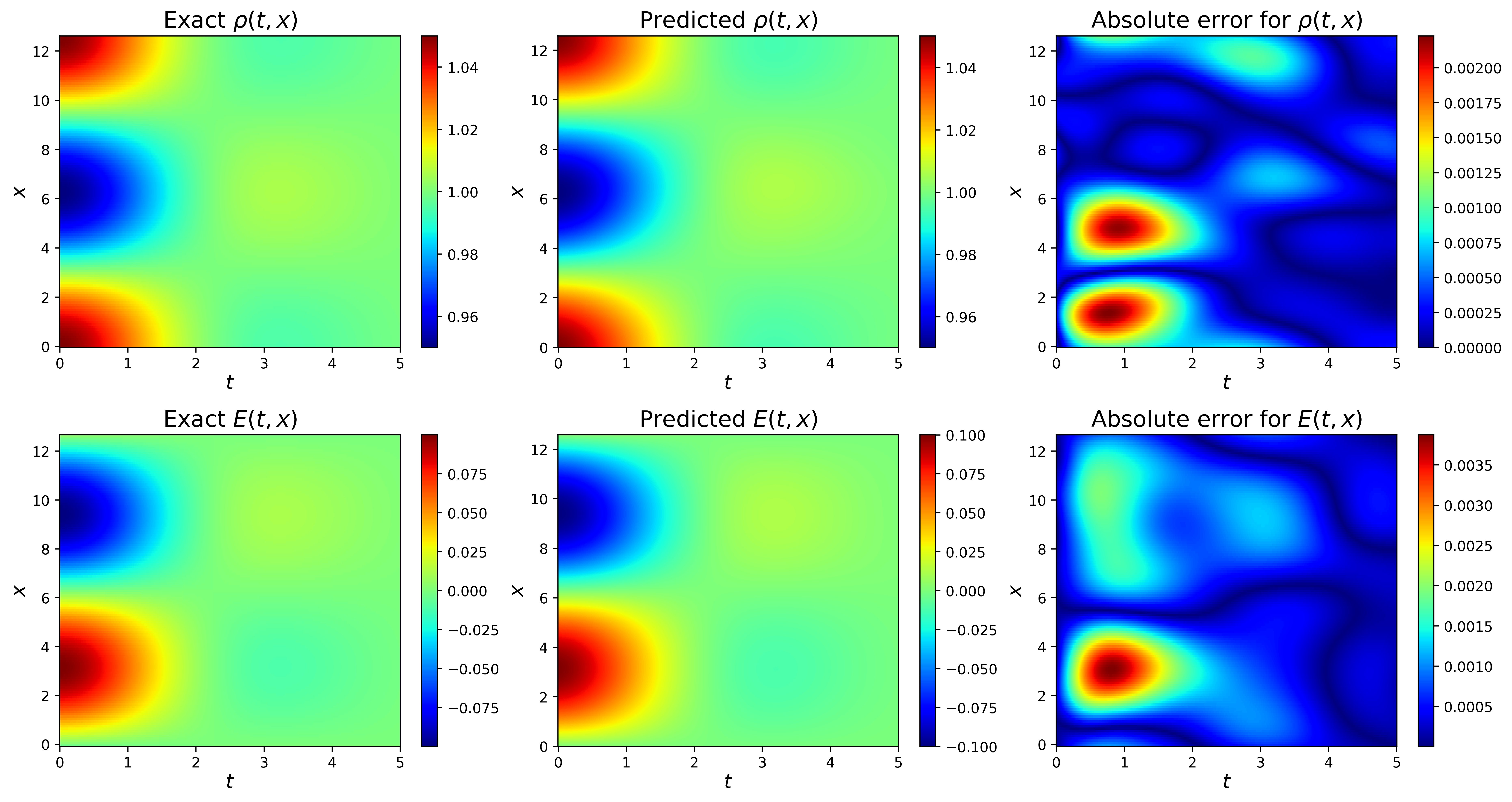}}
 
	\subfigure[Density $\rho$ (top) and electric field $E$ (bottom) in the high-field  regime $(\varepsilon=10^{-3})$.]{
	\includegraphics[width=0.7\linewidth]{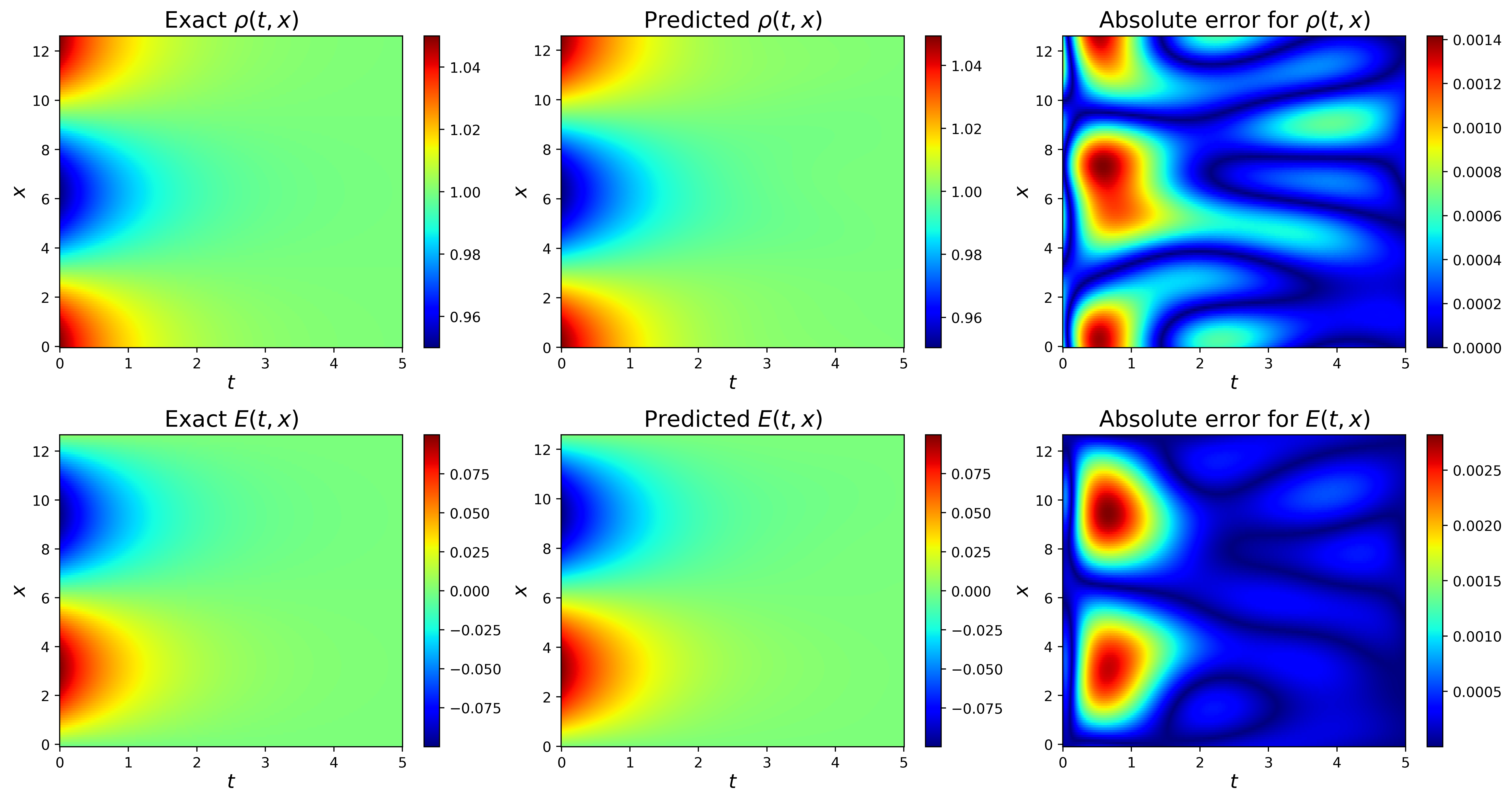}}
 
	\subfigure[Electric energy in the kinetic regime (left) and the high-field regime (right)]{
	\includegraphics[width=0.65\linewidth]{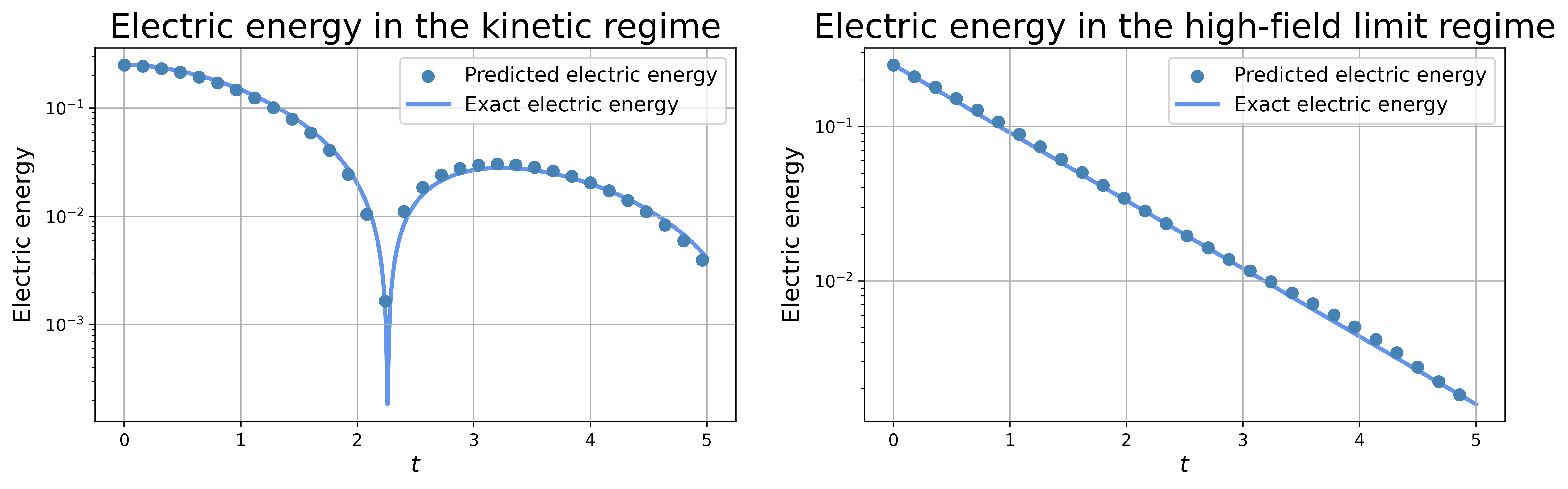}}
 
    \caption{Landau damping for the VPFP system solved by the AP-MIONet method in the kinetic and high-field regimes. Density $\rho$ and electric field $E$ plotted in $(t,x)\in[0,5]\times[0, 2\pi/k]$, and electric energy plotted in $t\in[0,5]$ for a representative input function. Penalty $\lambda_1=100$ and other penalties are set to $1$ for kinetic regime; Penalty $\lambda_1=500$ and other penalties are set to $1$ for the high-field regime.}
    \label{LD}
\end{figure}

We present the performance of the AP-MIONet method across multiscale regimes, characterized by the parameter $\varepsilon$.  We plot the density $\rho$ and the electric field $E$ as functions of time $t$ and space $x$, alongside the time evolution of the electric energy $\|E(t)\|_{L^2}$. Fig. \ref{LD} displays excellent agreement between the AP-MIONet predictions and the ground truth in both the kinetic and high-field regimes for a representative input sample from the test dataset. These results validate the AP property of the proposed AP-MIONet method, thereby supporting our theoretical claim. Moreover, we report the relative $\ell^2$-errors of density $\rho$, electric field $E$ and electric energy $\|E(t)\|_{L^2}$ over the entire spatio-temporal domain for the test dataset. The results for both the AP-MIONet and PI-MIONet methods are detailed in Table \ref{LD_error}. In this table, a dash ``-" indicates the non-convergence of the training process after reaching the designated training epochs. It is noteworthy that the relative $\ell^2$-errors for the PI-MIONet method are reported only up to $t=1.0$ instead of $t=5.0$ due to its deteriorating performance beyond this threshold. This table reveals that the PI-MIONet method exhibits limited effectiveness in the kinetic regime and fails to capture the asymptotic limit regime. We conclude that the AP-MIONet method significantly enhances the PI-MIONet method and proves to be an effective tool for accurately simulating multiscale plasma systems.

\vspace{10pt}
\begin{minipage}[!h]{\textwidth}
	\centering
	\begin{tabular}{cccc}
		\toprule
		$\ell^2$-error & $\rho$ & $E$ & Energy\\
		\midrule
		AP-MIONet & $0.11\%$  & $6.63\%$  & $7.06\%$\\
		PI-MIONet & $0.40\%$ & $9.81\%$ & $3.57\%$ \\
		\bottomrule 
        \multicolumn{4}{c}{(a) $\ell^2$-errors in the kinetic regime} 
	\end{tabular}
	\quad
	\begin{tabular}{cccc}
		\toprule
		$\ell^2$-error & $\rho$ & $E$ & Energy\\
		\midrule
		AP-MIONet & $0.08\%$  & $6.54\%$  & $9.94\%$\\
		PI-MIONet & $-$ & $-$ & $-$ \\
		\bottomrule 
        \multicolumn{4}{c}{(b) $\ell^2$-errors in the high-field regime}
	\end{tabular}\\
	\vspace{0.02\textwidth}
    \begin{minipage}[!h]{\textwidth}
        \centering
        Table $1$: For the Landau damping of the VPFP system, the $\ell^2$-errors of density $\rho$, electric field $E$ and electric energy $\|E(t)\|_{L^2}$ over the entire spatio-temporal domain for the test dataset using the AP-MIONet and PI-MIONet methods. The left table (a) shows the $\ell^2$-errors in the kinetic regime $(\varepsilon=1)$. The $\ell^2$-errors for the AP-MIONet method are reported up to $t=5.0$, while those for the PI-MIONet method are reported up to $t=1.0$. The right table (b) shows the $\ell^2$-errors in the high-field  regime $(\varepsilon=0.001)$. The dash ``-" indicates the non-convergence of the training process after reaching the designated training epochs.
    \end{minipage}
    \vspace{-0.005\textwidth}
\end{minipage}\label{LD_error}
\vspace{10pt}

\begin{remark}
    Temporal predictive capability of DeepONet (or MIONet)-based methods is one of the concerning problems. To extend predictions for a longer time, more sensors are necessary \cite{lu2021learning}, thereby demanding increased memory allocation. Recent studies \cite{lu2022solving, wang2023long} have introduced several enhancements to DeepONets to improve long-term predictive accuracy. These strategies can also be adapted to enhance the performance of the AP-MIONet method in long-time integrations.
\end{remark}

\subsection{Results for the semiconductor Boltzmann equations}

\subsubsection{The non-degenerate isotropic case}

In this case, the system is at low electron densities. When considering only the interaction of electrons with background impurities, the collision operator can be approximated by a linear relaxation time operator \cite{cercignani1997high}:
\begin{equation}
    \mathscr{Q}_{\text{isotropic}}(f)=\int M f^{\prime}-M^{\prime} f d v^{\prime}=M \rho-f,
\end{equation}
which represents the simplest case where the cross-section $\Psi\left(v^{\prime}, v\right)=1$ in Eq.\eqref{nond2}.  This model is commonly referred to as the ``time-relaxation" model.

\begin{figure}[!htb]
	\centering 
	\subfigbottomskip=5pt 
	\subfigcapskip=-5pt
	\subfigure[Density $\rho$ (top) and electric field $E$ (bottom) in the kinetic regime $(\varepsilon=1)$.]{
	\includegraphics[width=0.7\linewidth]{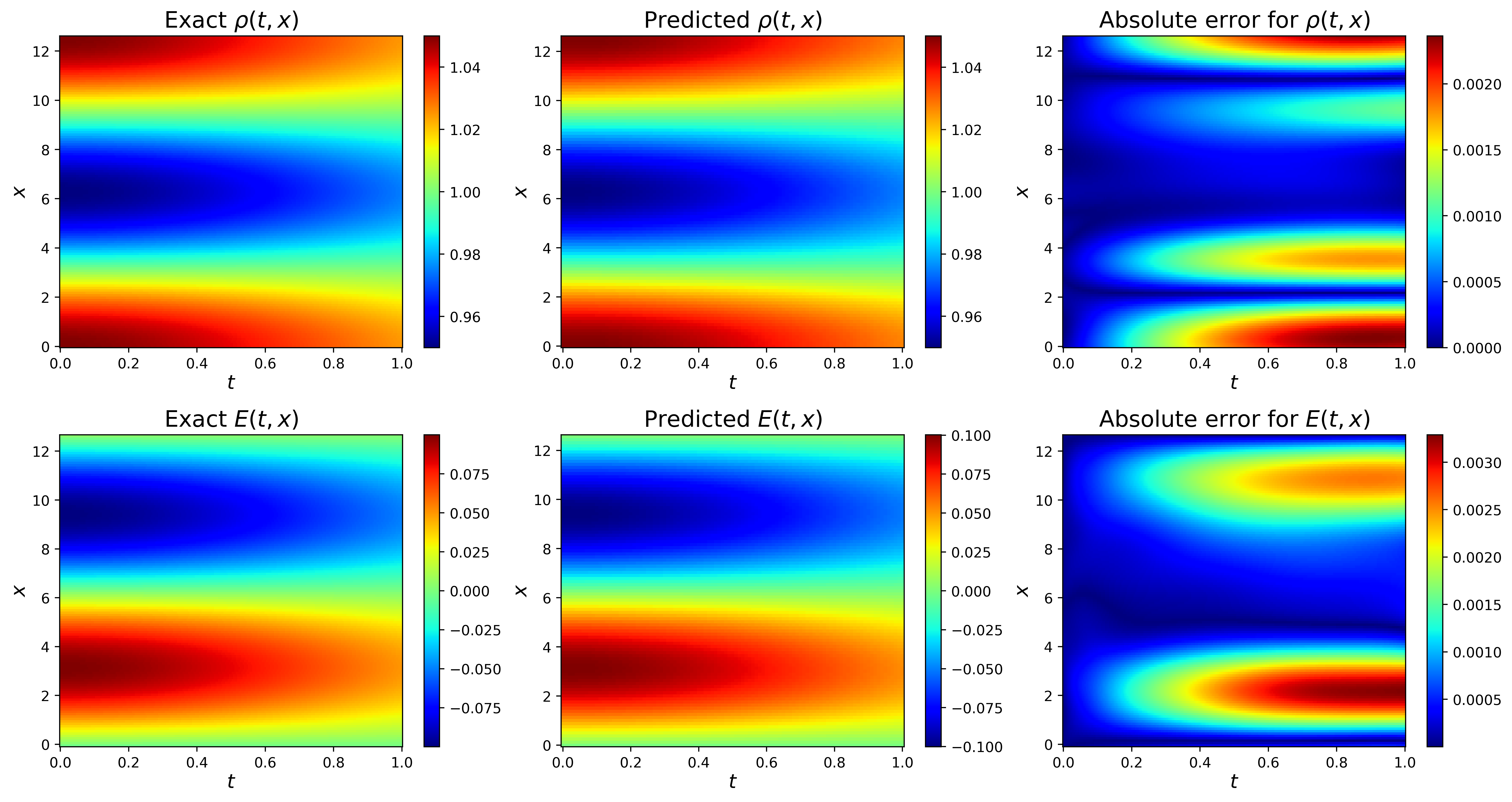}}
 
	\subfigure[Density $\rho$ (top) and electric field $E$ (bottom) in the high-field  regime $(\varepsilon=10^{-3})$.]{
	\includegraphics[width=0.7\linewidth]{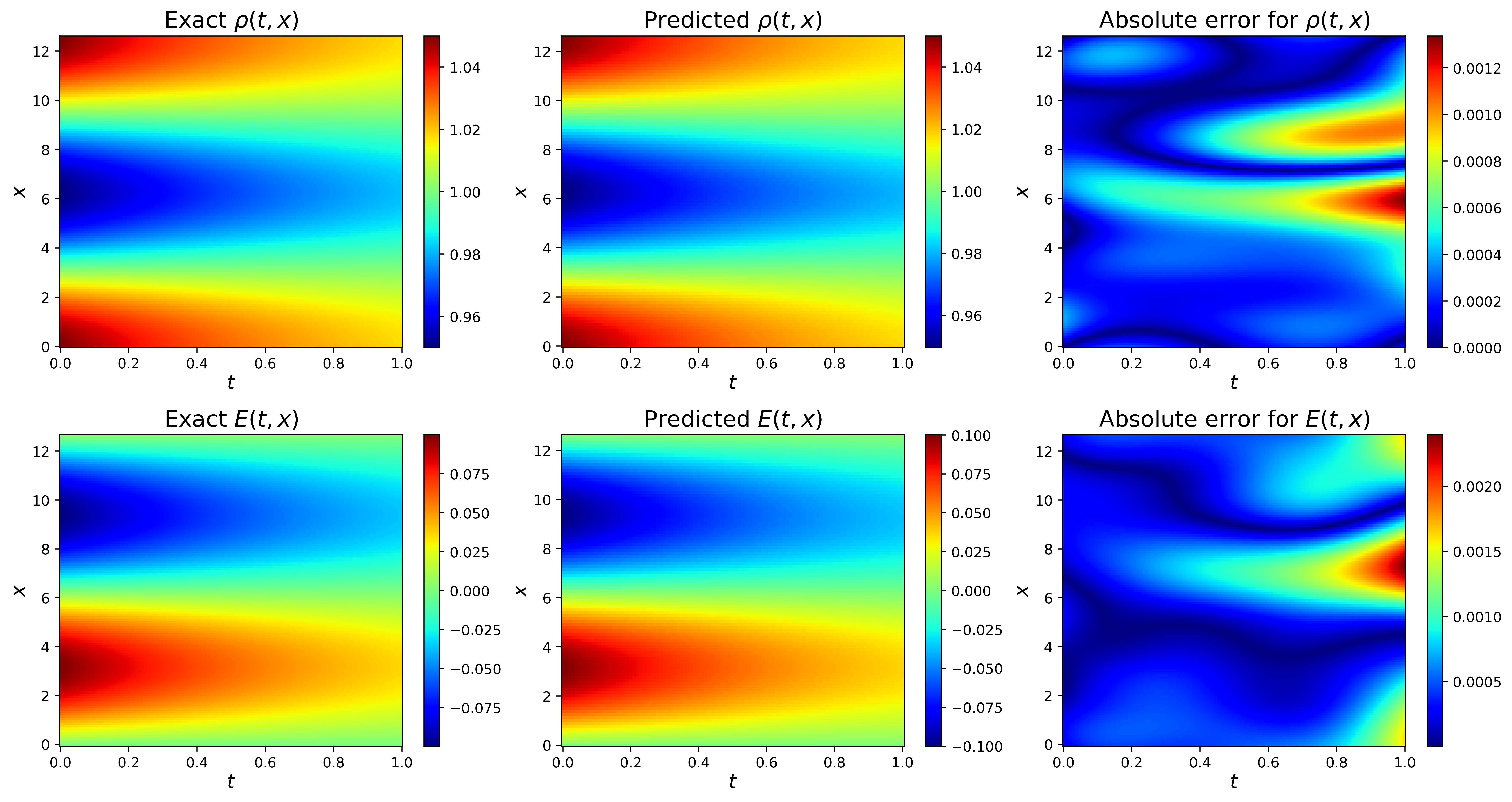}}
 
	\caption{Double peak instability for the non-degenerate isotropic semiconductor Boltzmann equation solved by the AP-MIONet method in the kinetic and high-field regimes. Density $\rho$ and electric field $E$ plotted in $(t,x)\in[0,1]\times[0, 2\pi/k]$ for a representative input function. Penalty $\lambda_1=20$ and other penalties are set to $1$ for the kinetic regime; Penalty $\lambda_1=100$ and other penalties are set to $1$ for the high-field regime. }
\label{DP}
\end{figure}

\vspace{1em}\noindent $\bullet$ \textbf{Double peak instability} \hspace{0.5em} We take this example to validate the capability of the AP-MIONet in accurately capturing the desired limiting macroscopic behavior for the non-degenerate isotropic case. The performance of the AP-MIONet in non-equilibrium states is also tested. The initial distribution is given by the sum of two local Maxwellian distributions:
\begin{equation}
    f_0(x, v)=\frac{h}{\sqrt{2 \pi}}\left( \frac{1}{2}\exp (-\frac{|v-1.5|^2}{2}) + \frac{1}{2}\exp (-\frac{|v+1.5|^2}{2})\right)(1+\alpha\cos(kx)) , \quad(x, v) \in[0,2 \pi / k] \times \mathbb{R},
\end{equation}
where $k$ is set to $0.5$ and perturbation parameter $\alpha$ is sampled uniformly from $[0.04, 0.06]$. The background charge $h$ is assumed to be uniform, with values sampled from the interval $[0.9, 1.1]$. Other settings remain consistent with those applied in the Landau damping case. Visualization results for a random input sample in test dataset are presented in Fig. \ref{DP}. The AP-MIONet predictions demonstrate strong concordance with the reference solutions. As shown in Fig. \ref{DP}(a), the AP-MIONet effectively adapts to non-equilibrium states. The relative $\ell^2$-errors of $\rho$ and $E$ in the kinetic regime are $0.12\%$ and $3.60\%$ respectively for the test dataset. In the high-field regime, these errors are $0.08\%$ for $\rho$ and $3.15\%$ for $E$. Above numerical results affirm that the AP-MIONet method effectively maintains the AP property for the time-relaxation model. 

\subsubsection{The non-degenerate anisotropic case}

In this part, we examine the non-degenerate anisotropic case with the collision cross-section defined as:
\begin{equation}
    \Psi\left(v, v^{\prime}\right)=1+e^{-\left(v-v^{\prime}\right)^2}.
\end{equation}
The collision operator $\mathscr{Q}_{\text{nond}}$ is then formulated by Eq.\eqref{nond1}.

\begin{figure}[!htb]
	\centering  
	\subfigbottomskip=5pt 
	\subfigcapskip=-5pt 
	\subfigure[Density $\rho$ (top) and electric field $E$ (bottom) in the kinetic regime $(\varepsilon=1)$.]{
	\includegraphics[width=0.7\linewidth]{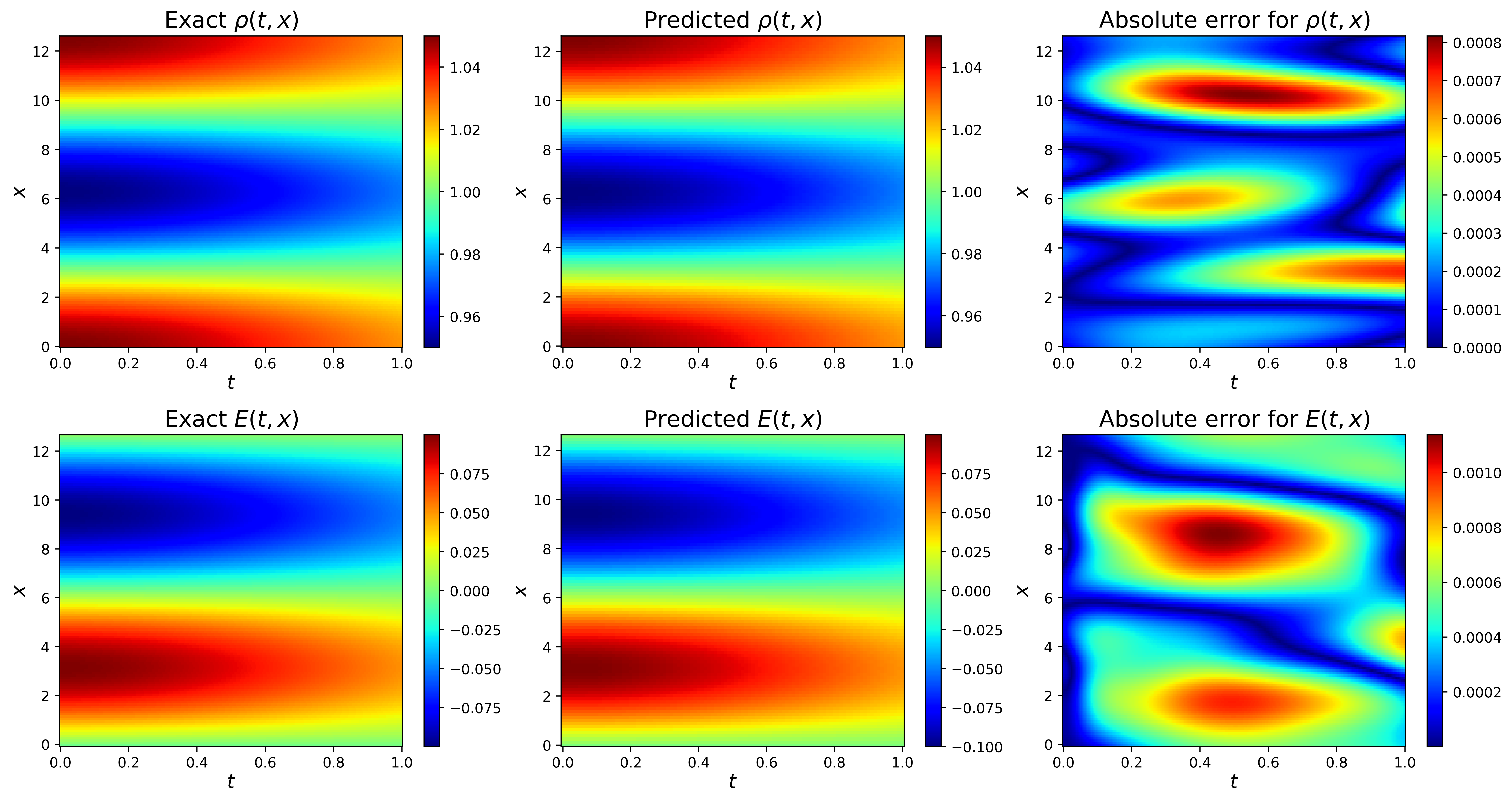}}
 
	\subfigure[Density $\rho$ (top) and electric field $E$ (bottom) in the high-field regime $(\varepsilon=10^{-3})$.]{
	\includegraphics[width=0.7\linewidth]{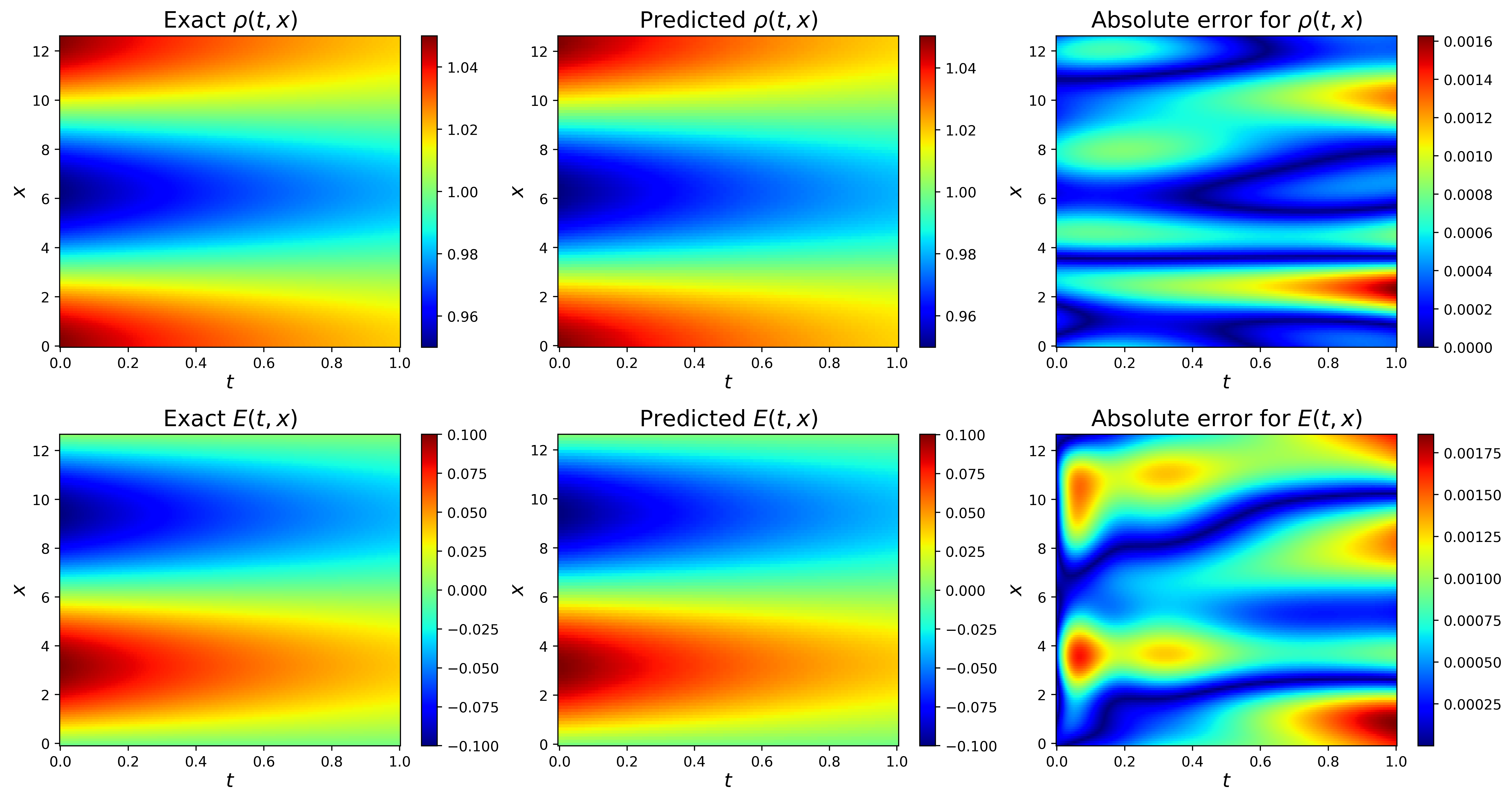}}
 
	\caption{Two stream instability for the non-degenerate anisotropic semiconductor Boltzmann equation solved by the AP-MIONet method in the kinetic and high-field regimes. Density $\rho$ and electric field $E$ plotted in $(t,x)\in[0,1]\times[0, 2\pi/k]$ for a representative input function. Penalty $\lambda_1=20$ and other penalties are set to $1$ for the kinetic regime; Penalty $\lambda_1=100$ and other penalties are set to $1$ for the high-field regime.}
\label{TS}
\end{figure}

\vspace{1em}\noindent $\bullet$ \textbf{Two stream instability} \hspace{0.5em} We utilize this instance to evaluate the performance of the AP-MIONet to accurately capture  the solution to the high-field limit equation. The initial distribution takes the following form:
\begin{equation}
    f_0(x, v)=\frac{h}{\sqrt{2 \pi}} v^2 \exp \left(-\frac{v^2}{2}\right)(1+\alpha \cos (k x)),
\end{equation}
where $k$ is set to $0.5$ and $\alpha$ is sampled from $[0.04, 0.06]$. The background charge $h$ is uniform and sampled from $[0.9,1.1]$. Other settings are the same as those applied in the Landau damping case. Fig. \ref{TS} exhibits a comparison between the predicted and reference solutions for a random input sample from the test dataset. We observe significant alignment between the AP-MIONet predictions and the reference solutions in both the kinetic and high-field regimes. The relative $\ell^2$-errors of $\rho$ and $E$ for the test dataset are $0.10\%$ and $3.62\%$ in the kinetic regime, respectively. In the high-field regime, these errors are $0.12\%$ for $\rho$ and $4.82\%$ for $E$. These findings validate the AP property of the AP-MIONet method in the non-degenerate anisotropic case. 

\subsubsection{The degenerate case}

This part addresses the degenerate case, where the collision operator $\mathscr{Q}_{\text{deg}}$ is defined in Eq. \eqref{deg} with the cross-section $\Psi\left(v, v^{\prime}\right)=1+e^{-\left(v-v^{\prime}\right)^2}$. We note that the operator $\mathscr{Q}_{\text{deg}}$ is nonlinear, and the corresponding limit equation is also a nonlinear convection equation. Such nonlinearity may pose a challenge for neural network-based methods in accurately solving the governing PDEs \cite{li2022model, wang2021learning}. In this part, we will investigate the ability of the AP-MIONet method to handle the nonlinear challenge.

\begin{figure}[!htb]
	\centering 
	\subfigbottomskip=5pt
	\subfigcapskip=-5pt 
	\subfigure[Density $\rho$ (top) and electric field $E$ (bottom) in the kinetic regime $(\varepsilon=1)$.]{
	\includegraphics[width=0.7\linewidth]{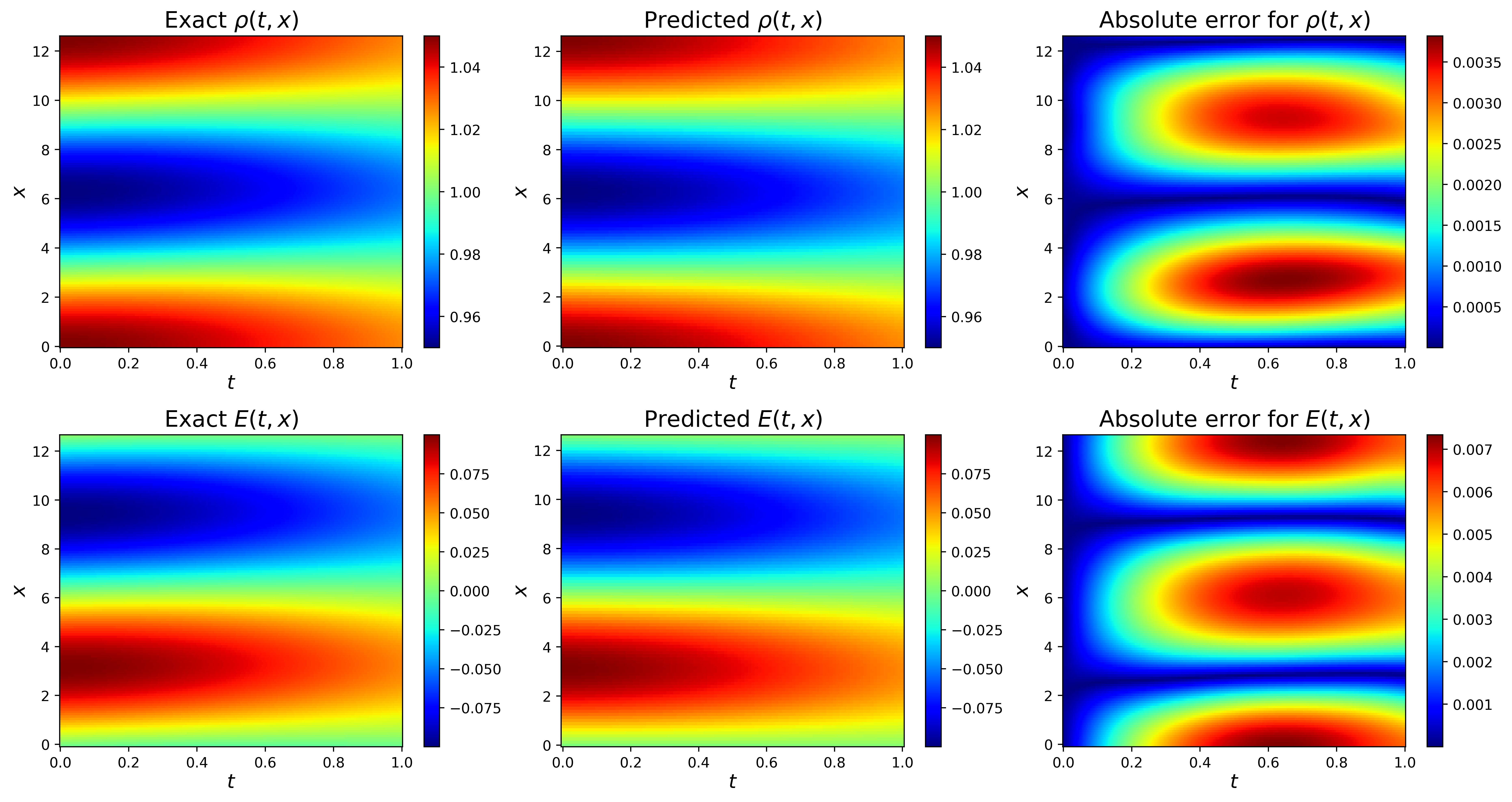}}
 
	\subfigure[Density $\rho$ (top) and electric field $E$ (bottom) in the high-field regime $(\varepsilon=10^{-3})$.]{
	\includegraphics[width=0.7\linewidth]{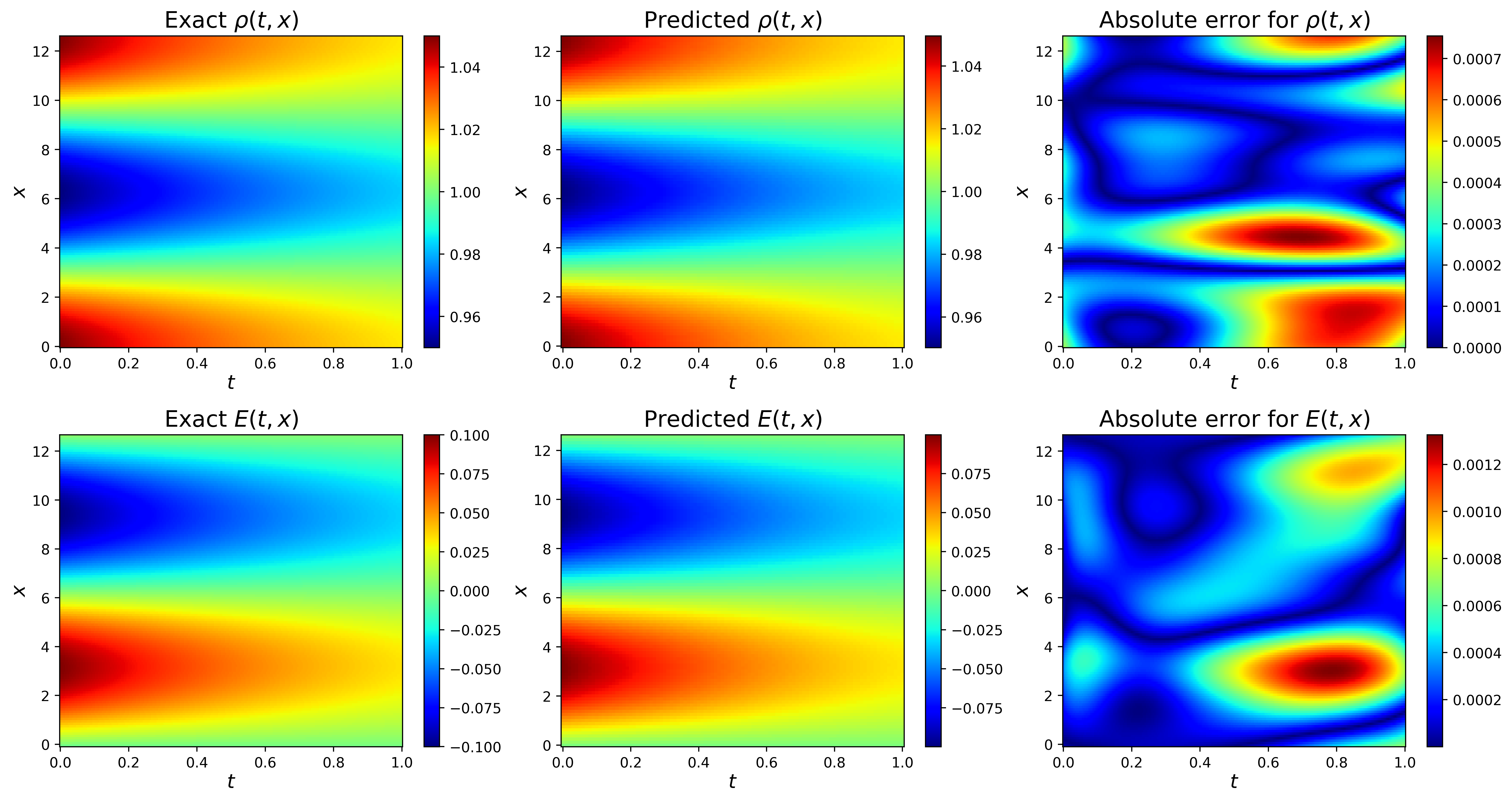}}
 
	\caption{Bump-on-tail instability for the degenerate semiconductor Boltzmann equation solved by the AP-MIONet method in the kinetic and high-field regimes. Density $\rho$ and electric field $E$ plotted in $(t,x)\in[0,1]\times[0, 2\pi/k]$ for a representative input function. Penalty $\lambda_1=20$ and other penalties are set to $1$ for the kinetic regime; Penalty $\lambda_1=100$ and other penalties are set to $1$ for the high-field regime.}
\label{BOT}
\end{figure}

\vspace{1em}\noindent $\bullet$ \textbf{Bump-on-tail instability} \hspace{0.5em} This example serves to verify the AP property of the AP-MIONet method and its ability to manage nonlinearity in the degenerate case. The initial distribution is configured to represent a bump on the tail of the Maxwellian distribution:
\begin{equation}
    f_0(x, v)=\frac{h}{\sqrt{2 \pi}}\left(\frac{9}{10} \exp (-\frac{v^2}{2})+\frac{2}{10} \exp \left(-4(v-4.5)^2\right)\right)(1+\alpha \cos (k x)),
\end{equation}
where $k$ is set to $0.5$ and $\alpha$ is sampled from $[0.04,0.06]$. The background charge $h$ is uniform and sampled from $[0.9,1.1]$. The computational domain is $[0,2 \pi / k] \times\left[v_{\min }, v_{\max }\right]$, where $-v_{\min }=v_{\max }=8$. Other settings remain the same as those used in the Landau damping case. As shown in Fig. \ref{BOT}, remarkable consistency is achieved between the predicted and reference solutions for a representative example in the test dataset. The relative $\ell^2$-errors of $\rho$ and $E$ for the test dataset are $0.23\%$ and $8.00\%$ in the kinetic regime, respectively. In the high-field regime, these errors are $0.06\%$ for $\rho$ and $2.61\%$ for $E$. Therefore, the AP property of our method can be confirmed by the visualizations and predictive errors, along with its ability in tackling nonlinearity. 

\vspace{1em}\noindent $\bullet$ \textbf{Mixing regimes} \hspace{0.5em} The AP-MIONet method is now evaluated in the mixing regimes, where the scale parameter $\varepsilon$ exhibits spatial variations across several orders of magnitude. The parameter $\varepsilon$ is formulated as:
\begin{equation}
    \varepsilon(x)= \begin{cases}\varepsilon_0+\frac{1}{2}(\tanh (5-10 x)+\tanh (5+10 x)), & -1 \leq x \leq 0.3, \\ \varepsilon_0, & 0.3<x \leq 1,\end{cases}
\end{equation}
where $\varepsilon_0$ is set to $0.001$. This setting considers both the kinetic and high-field regimes. The initial distribution is assumed to be:
\begin{equation}
    f_0(x, v)=\frac{h}{2\sqrt{2 \pi}} \exp \left(-\frac{v^2}{2}\right)(2+\sin(k x)),
\end{equation}
where $k$ is set as $0.5$ and the uniform $h$ is sampled from the range $[0.80,0.85]$. The computational domain is $[-1, 1] \times\left[v_{\min }, v_{\max }\right]$ with $-v_{\min }=v_{\max }=6$. The periodic BC in $x$-direction is imposed. Our goal is to learn three solution operators $\mathcal{P}, \mathcal{F}$ and $\Phi$ given by:
\begin{equation}
    \mathcal{F}:h \mapsto f, \quad\mathcal{P}:h \mapsto\rho, \quad\Phi: h \mapsto\phi.
\end{equation}
To achieve this, each MIONet degenerates to include one branch net and one trunk net. Additionally, we adopt the early stopping strategy \cite{prechelt2002early, yao2007early} to avoid overfitting during the training process.

\begin{figure}[!htb]
	\centering 
	\subfigbottomskip=5pt
	\subfigcapskip=-5pt 
	\subfigure[The scale parameter $\varepsilon$ as a function of space $x$ with several orders of magnitude]{
	\includegraphics[width=0.55\linewidth]{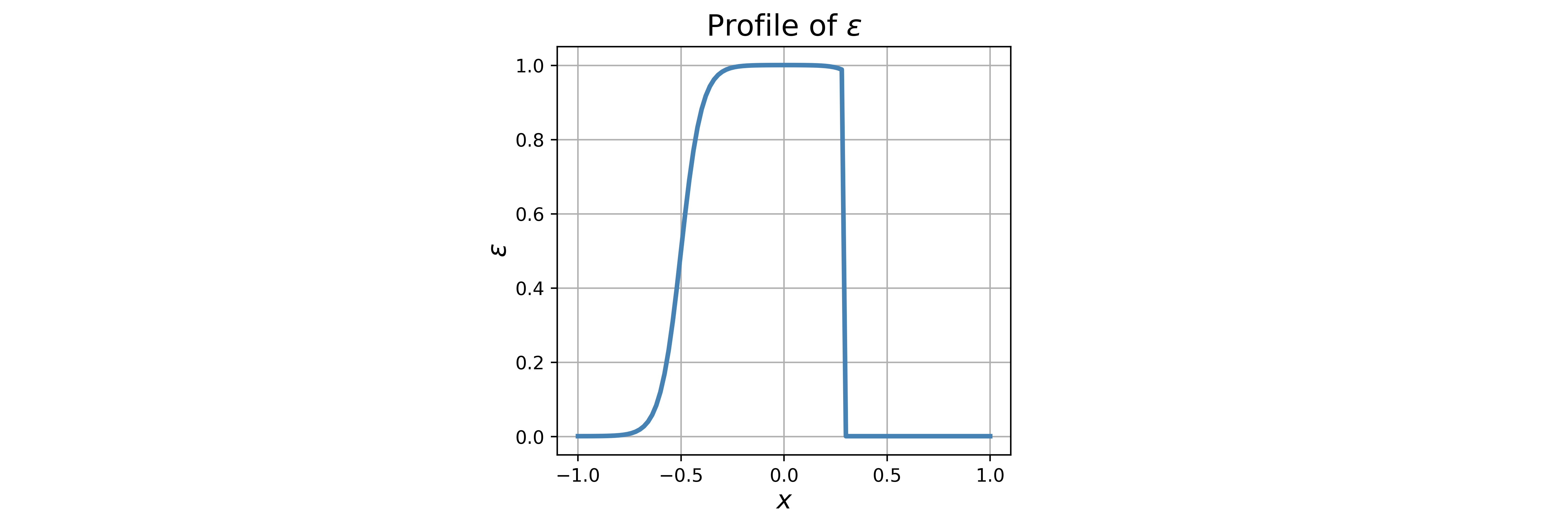}}
 
	\subfigure[Density $\rho$ (top) and electric field $E$ (bottom) in the mixing regimes.]{
	\includegraphics[width=0.7\linewidth]{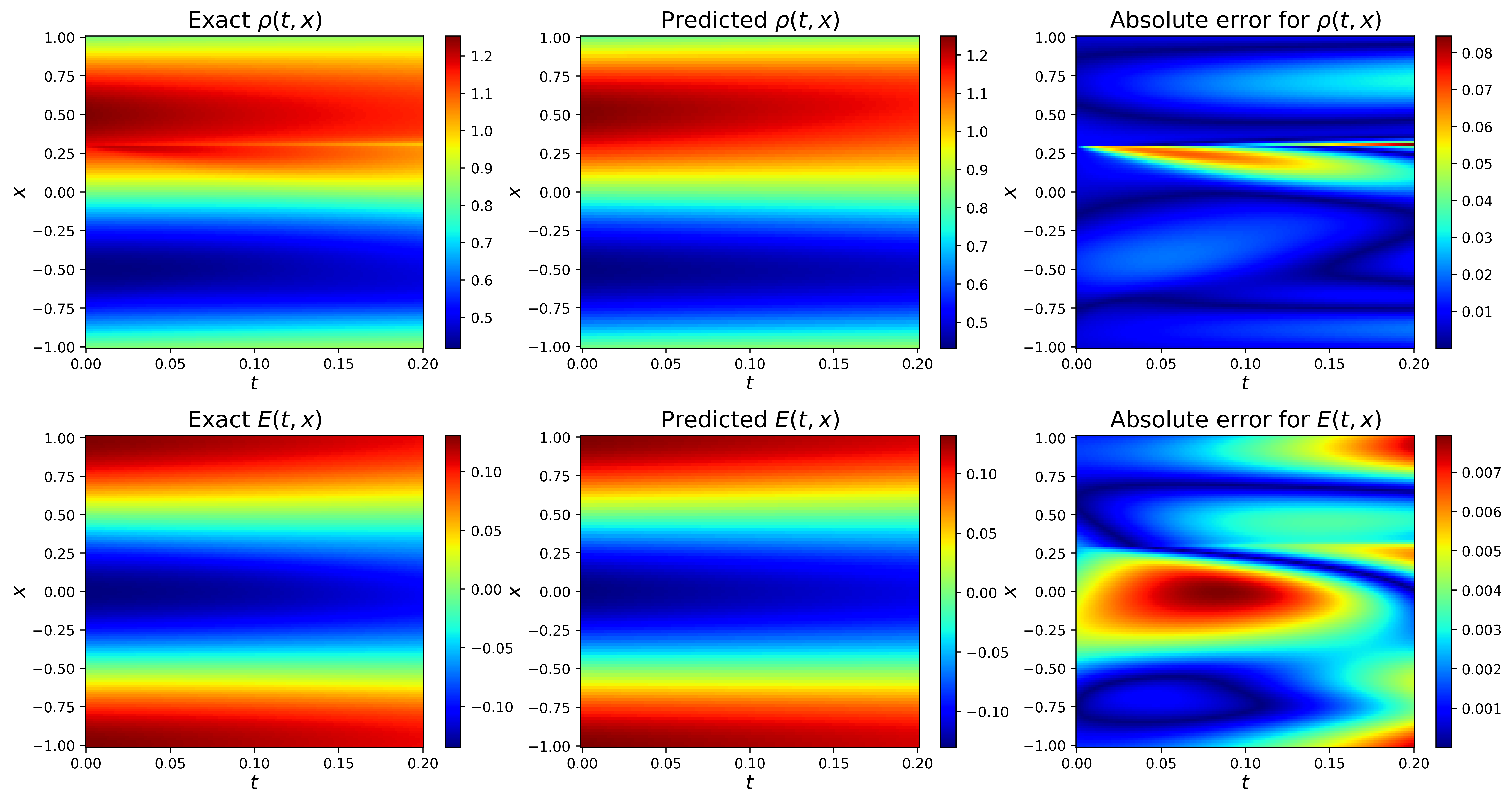}}

    \subfigure[Cross sections of $\rho$ and $E$ at $t=0.0, 0.1, 0.2$.]{
	\includegraphics[width=0.65\linewidth]{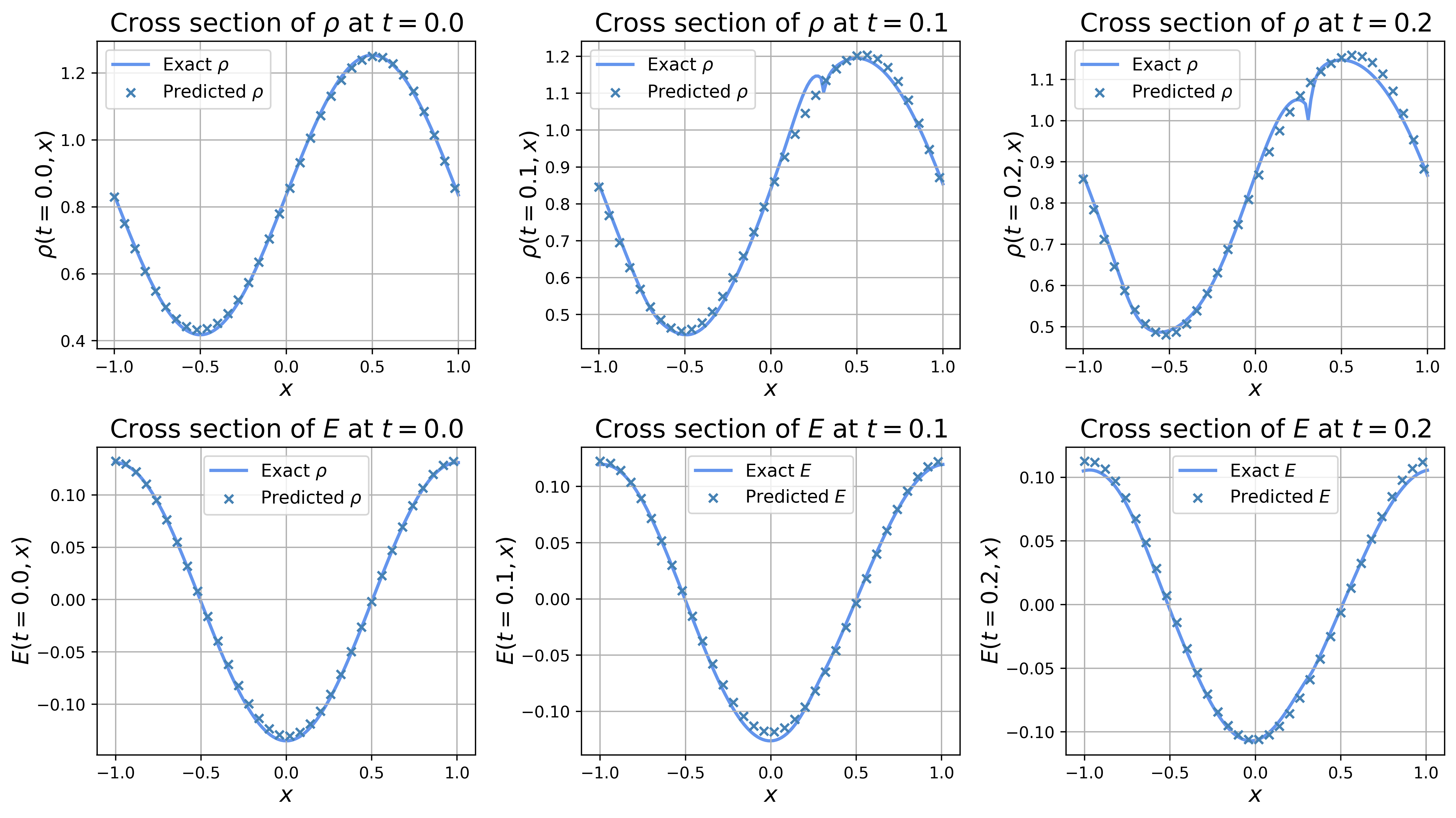}}
 
	\caption{Mixing regimes problem for the degenerate semiconductor Boltzmann equation solved by the AP-MIONet method. The subfigure (a) shows the profile of parameter $\varepsilon$ with the mixing regimes. In subfigure (b), density $\rho$ and electric field $E$ are plotted in $(t,x)\in[0,0.2]\times[-1, 1]$ for a representative input function. All penalties  are set to $1$. The subfigure (c) displays the cross sections of $\rho$ and $E$ with $t=0.0, 0.1, 0.2$ for visualization purposes. }
    \label{MIX}
\end{figure}

We assess the capability of the AP-MIONet to tackle the case involving mixing regimes. As shown in panel (a) of Fig. \ref{MIX}, the scale parameter $\varepsilon$ exhibits a discontinuity at $x=0.3$. Panels (b) and (c) of Fig. \ref{MIX} illustrate that the AP-MIONet captures the profile of the density $\rho$ with notable accuracy except for a single point of minor discontinuity at $x=0.3$. The inability in capturing the localized discontinuity can be attributed to the inherent difficulty of neural networks in learning high-frequency information, a challenge arising from the spectral bias \cite{rahaman2019spectral}. The relative $\ell^2$-errors of $\rho$ and $E$ are $2.08\%$ and $4.15\%$ respectively for the test dataset. 

\subsection{Computational cost}

\begin{figure}[!htb]
	\centering 
	\includegraphics[width=0.40\linewidth]{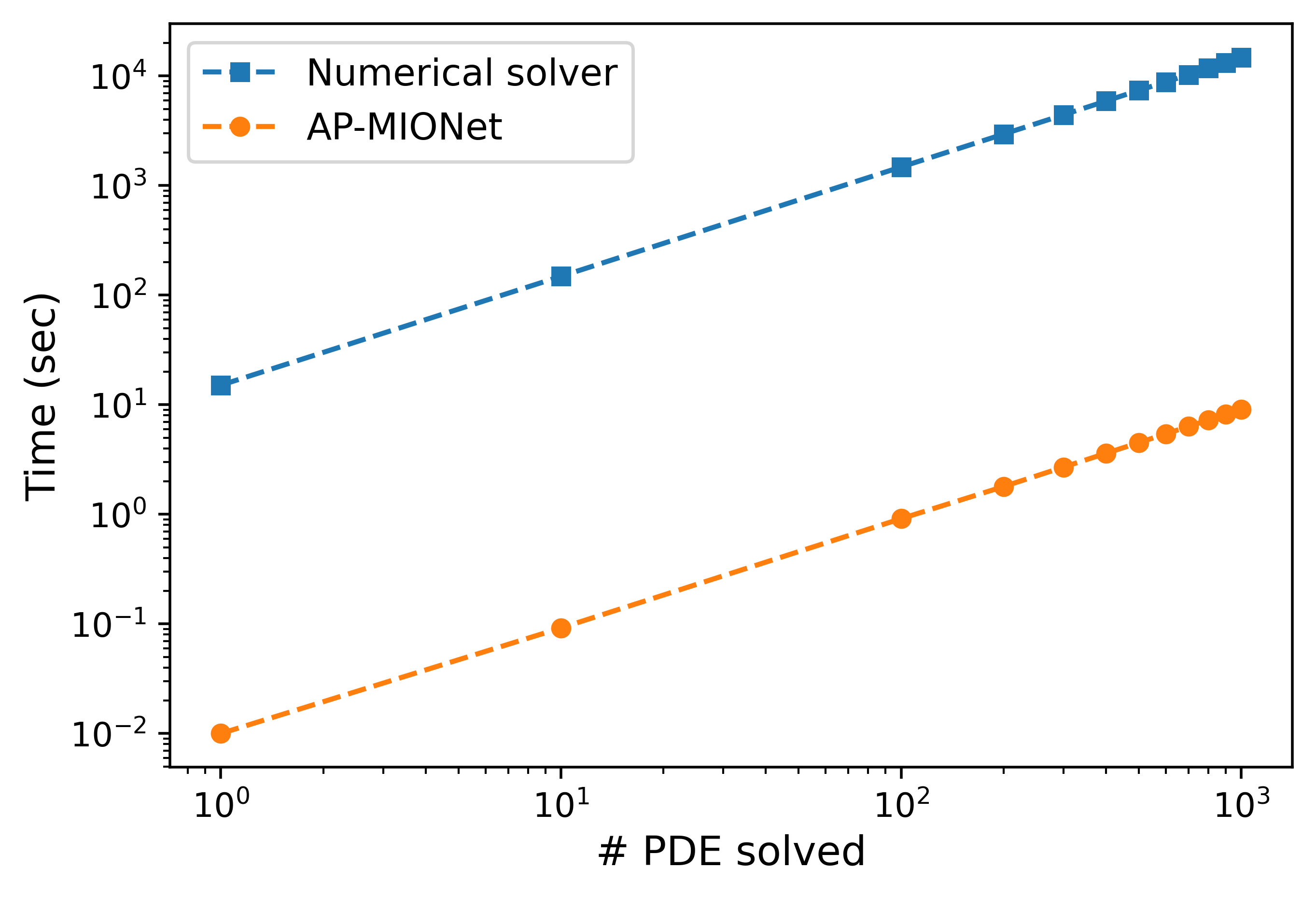}
	\caption{Computational cost for performing inference for the bump-on-tail instability example with a trained AP-MIONet, compared to a conventional asymptotic-preserving solver \cite{jin2013asymptotic}. Reported timings are obtained on a single NVIDIA TITAN Xp GPU.}
\label{infer}
\end{figure}

This section shows  the computational costs for the AP-MIONet. Table \ref{Cost} in Appendix B summarizes the hours spent on training the PI-MIONet and AP-MIONet models for the examples in Sect. 5. All networks are trained in parallel using four NVIDIA TITAN Xp graphics processing units (GPUs). We then compare the computational costs of inference between the AP-MIONet and a conventional numerical solver. As an illustrative example, we analyze the bump-on-tail instability for solving the degenerate semiconductor Boltzmann equations in the high-field regime ($\varepsilon=0.001$). The AP-MIONet model is implemented using the JAX library. For comparison, a traditional AP numerical scheme \cite{jin2013asymptotic} is executed on one GPU using the CuPy library. As depicted in Fig.\ref{infer}, a trained AP-MIONet can predict the solutions of $\mathcal{O}\left(10^3\right)$ PDEs in one second. This remarkable performance yields a speed-up of three orders of magnitude compared to the traditional AP numerical solver. Moreover, inference with AP-MIONets is trivially parallelizable. This parallelization strategy can leverage multiple GPUs to simultaneously predict solutions for an expanded set of PDEs, thereby further enhancing the computational efficiency.

This is the first study for such problems and we will next study it for multi-dimensional problems. 

\section{Conclusion}
In this paper, we focus on learning the solution operators for multiscale kinetic equations with the high-field scaling, specifically the Vlasov-Poisson-Fokker-Planck system and the semiconductor Boltzmann equation. To address the additional variable parameter derived from the coupled electric field, we employ the multiple-input DeepONets (MIONets) to parameterize the solution operators. Facing up the multiscale challenge, the Asymptotic-preserving (AP) mechanism is integrated into the loss function of the MIONets. We reformulate the loss function by incorporating both the mass conservation law and the original kinetic equation, enabling it to preserve mass and adapt to non-equilibrium states. This AP loss does not require an explicit local equilibrium, promising its feasibility for the semiconductor Boltzmann equation. Therefore, we propose a new Asymptotic-preserving MIONet (AP-MIONet) method, a generic framework for solving these two kinetic equations with the high-field scaling. Extensive numerical results validate the effectiveness of our approach. The AP-MIONet method extends the AP concept within operator learning for a realistic and complex scenario where the external fields are prevalent and markedly influence particle dynamics.

\section*{Acknowledgement}

This research was supported  by the Strategic Priority Research Program of Chinese
Academy of Sciences XDA25010401.

\section*{Data availability statement}
Data availability is not applicable to this article as no new data were generated or analyzed in this research.

\section*{Declarations}
\textbf{Conflict of interest:} No potential conflict of interest is declared by the authors.

\bibliographystyle{siam}
\bibliography{ref.bib}

\begin{thebibliography}{10}

\bibitem{abdallah2001high}
{\sc N.~B. Abdallah and H.~Chaker}, {\em The high field asymptotics for degenerate semiconductors}, Mathematical Models and Methods in Applied Sciences, 11 (2001), pp.~1253--1272.

\bibitem{abdallah2000convergence}
{\sc N.~B. Abdallah, L.~Desvillettes, and S.~G{\'e}nieys}, {\em On the convergence of the {B}oltzmann equation for semiconductors toward the energy transport model}, Journal of Statistical Physics, 98 (2000), pp.~835--870.

\bibitem{li2020neural}
{\sc A.~Anandkumar, K.~Azizzadenesheli, K.~Bhattacharya, N.~Kovachki, Z.~Li, B.~Liu, and A.~Stuart}, {\em {N}eural {O}perator: {G}raph {K}ernel {N}etwork for {P}artial {D}ifferential {E}quations}, in ICLR 2020 Workshop on Integration of Deep Neural Models and Differential Equations, 2019.

\bibitem{paszke2017automatic}
{\sc A.~G. Baydin, B.~A. Pearlmutter, A.~A. Radul, and J.~M. Siskind}, {\em Automatic differentiation in machine learning: a survey}, Journal of Machine Learning Research, 18 (2018), pp.~1--43.

\bibitem{ben2007high}
{\sc N.~Ben~Abdallah, H.~Chaker, and C.~Schmeiser}, {\em The high field asymptotics for a fermionic {B}oltzmann equation: entropy solutions and kinetic shock profiles}, Journal of Hyperbolic Differential Equations, 4 (2007), pp.~679--704.

\bibitem{ben1996energy}
{\sc N.~Ben~Abdallah, P.~Degond, and S.~G{\'e}nieys}, {\em An energy-transport model for semiconductors derived from the {B}oltzmann equation}, Journal of statistical physics, 84 (1996), pp.~205--231.

\bibitem{bertaglia2022asymptotic}
{\sc G.~Bertaglia, C.~Lu, L.~Pareschi, and X.~Zhu}, {\em Asymptotic-{P}reserving {N}eural {N}etworks for multiscale hyperbolic models of epidemic spread}, Mathematical Models and Methods in Applied Sciences, 32 (2022), pp.~1949--1985.

\bibitem{blakemore2002semiconductor}
{\sc J.~S. Blakemore}, {\em Semiconductor statistics}, Courier Corporation, 2002.

\bibitem{bouchut1993existence}
{\sc F.~Bouchut}, {\em Existence and uniqueness of a global smooth solution for the {V}lasov-{P}oisson-{F}okker-{P}lanck system in three dimensions}, Journal of functional analysis, 111 (1993), pp.~239--258.

\bibitem{cercignani1997high}
{\sc C.~Cercignani, I.~M. Gamba, and C.~D. Levermore}, {\em High field approximations to a {B}oltzmann-{P}oisson system and boundary conditions in a semiconductor}, Applied Mathematics Letters, 10 (1997), pp.~111--117.

\bibitem{chandrasekhar1943stochastic}
{\sc S.~Chandrasekhar}, {\em Stochastic problems in physics and astronomy}, Reviews of Modern Physics, 15 (1943), p.~1.

\bibitem{chen1995universal}
{\sc T.~Chen and H.~Chen}, {\em Universal approximation to nonlinear operators by neural networks with arbitrary activation functions and its application to dynamical systems}, IEEE transactions on neural networks, 6 (1995), pp.~911--917.

\bibitem{chen2022solving}
{\sc Z.~Chen, L.~Liu, and L.~Mu}, {\em Solving the linear transport equation by a deep neural network approach}, Discrete and Continuous Dynamical Systems-S, 15 (2022), pp.~669--686.

\bibitem{crouseilles2011asymptotic}
{\sc N.~Crouseilles and M.~Lemou}, {\em An asymptotic preserving scheme based on a micro-macro decomposition for collisional {V}lasov equations: diffusion and high-field scaling limits.}, Kinetic and related models, 4 (2011), pp.~441--477.

\bibitem{degend1986global}
{\sc P.~Degend}, {\em Global existence of smooth solution for {V}lasov-{F}okker-{P}lank equation in 1 and 2 space dimension}, Annales scientifiques de l'{\'E}cole normale sup{\'e}rieure, 19 (1986), pp.~519--542.

\bibitem{dong2021method}
{\sc S.~Dong and N.~Ni}, {\em A method for representing periodic functions and enforcing exactly periodic boundary conditions with deep neural networks}, Journal of Computational Physics, 435 (2021), p.~110242.

\bibitem{fanaskov2023spectral}
{\sc V.~S. Fanaskov and I.~V. Oseledets}, {\em Spectral neural operators}, in Doklady Mathematics, vol.~108, Springer, 2023, pp.~S226--S232.

\bibitem{frosali1989scattering}
{\sc G.~Frosali and C.~V.~M. van~der Mee}, {\em Scattering theory relevant to the linear transport of particle swarms}, Journal of statistical physics, 56 (1989), pp.~139--148.

\bibitem{frosali1989conditions}
{\sc G.~Frosali, C.~V.~M. van~der Mee, and S.~L. Paveri-Fontana}, {\em Conditions for runaway phenomena in the kinetic theory of particle swarms}, Journal of mathematical physics, 30 (1989), pp.~1177--1186.

\bibitem{glorot2010understanding}
{\sc X.~Glorot and Y.~Bengio}, {\em Understanding the difficulty of training deep feedforward neural networks}, in Proceedings of the thirteenth international conference on artificial intelligence and statistics, JMLR Workshop and Conference Proceedings, 2010, pp.~249--256.

\bibitem{golse1992limite}
{\sc F.~Golse and F.~Poupaud}, {\em Limite fluide des {\'e}quations de {B}oltzmann des semi-conducteurs pour une statistique de {F}ermi--{D}irac}, Asymptotic Analysis, 6 (1992), pp.~135--160.

\bibitem{goswami2023physics}
{\sc S.~Goswami, A.~Bora, Y.~Yu, and G.~E. Karniadakis}, {\em Physics-informed deep neural operator networks}, in Machine Learning in Modeling and Simulation: Methods and Applications, Springer, 2023, pp.~219--254.

\bibitem{goudon2005multidimensional}
{\sc T.~Goudon, J.~Nieto, F.~Poupaud, and J.~Soler}, {\em Multidimensional high-field limit of the electrostatic {V}lasov--{P}oisson--{F}okker--{P}lanck system}, Journal of Differential Equations, 213 (2005), pp.~418--442.

\bibitem{han2018solving}
{\sc J.~Han, A.~Jentzen, and W.~E}, {\em Solving high-dimensional partial differential equations using deep learning}, Proceedings of the National Academy of Sciences, 115 (2018), pp.~8505--8510.

\bibitem{hornik1990universal}
{\sc K.~Hornik, M.~Stinchcombe, and H.~White}, {\em Universal approximation of an unknown mapping and its derivatives using multilayer feedforward networks}, Neural networks, 3 (1990), pp.~551--560.

\bibitem{hu2024hybrid}
{\sc J.~Hu and P.~Jin}, {\em A hybrid iterative method based on {MION}et for {PDE}s: {T}heory and numerical examples}, arXiv preprint arXiv:2402.07156,  (2024).

\bibitem{hwang2020trend}
{\sc H.~J. Hwang, J.~W. Jang, H.~Jo, and J.~Y. Lee}, {\em Trend to equilibrium for the kinetic {F}okker-{P}lanck equation via the neural network approach}, Journal of Computational Physics, 419 (2020), p.~109665.

\bibitem{jiang2023fourier}
{\sc Z.~Jiang, M.~Zhu, D.~Li, Q.~Li, Y.~O. Yuan, and L.~Lu}, {\em Fourier-{MION}et: {F}ourier-enhanced multiple-input neural operators for multiphase modeling of geological carbon sequestration}, arXiv preprint arXiv:2303.04778,  (2023).

\bibitem{jin2022mionet}
{\sc P.~Jin, S.~Meng, and L.~Lu}, {\em {MION}et: Learning multiple-input operators via tensor product}, SIAM Journal on Scientific Computing, 44 (2022), pp.~A3490--A3514.

\bibitem{jin2022asymptotic}
{\sc S.~Jin}, {\em Asymptotic-preserving schemes for multiscale physical problems}, Acta Numerica, 31 (2022), pp.~415--489.

\bibitem{jin2023asymptotic1}
{\sc S.~Jin, Z.~Ma, and K.~Wu}, {\em Asymptotic-preserving neural networks for multiscale time-dependent linear transport equations}, Journal of Scientific Computing, 94 (2023), p.~57.

\bibitem{jin2023asymptotic}
\leavevmode\vrule height 2pt depth -1.6pt width 23pt, {\em Asymptotic-{P}reserving {N}eural {N}etworks for {M}ultiscale {K}inetic {E}quations}, Communications in Computational Physics, 35 (2024), pp.~693--723.

\bibitem{jin2023asymptotic2}
{\sc S.~Jin, Z.~Ma, and T.-a. Zhang}, {\em Asymptotic-{P}reserving {N}eural {N}etworks for {M}ultiscale {V}lasov--{P}oisson--{F}okker--{P}lanck {S}ystem in the {H}igh-{F}ield {R}egime}, Journal of Scientific Computing, 99 (2024), p.~61.

\bibitem{jin2011asymptotic}
{\sc S.~Jin and L.~Wang}, {\em An asymptotic preserving scheme for the {V}lasov-{P}oisson-{F}okker-{P}lanck system in the high field regime}, Acta Mathematica Scientia, 31 (2011), pp.~2219--2232.

\bibitem{jin2013asymptotic}
\leavevmode\vrule height 2pt depth -1.6pt width 23pt, {\em Asymptotic-preserving numerical schemes for the semiconductor {B}oltzmann equation efficient in the high field regime}, SIAM Journal on Scientific Computing, 35 (2013), pp.~B799--B819.

\bibitem{kingma2014adam}
{\sc D.~P. Kingma and J.~Ba}, {\em Adam: A method for stochastic optimization}, arXiv preprint arXiv:1412.6980,  (2014).

\bibitem{kong2023b}
{\sc Z.~Kong, A.~Mollaali, C.~Moya, N.~Lu, and G.~Lin}, {\em {B}-{LSTM}-{MION}et: {B}ayesian {LSTM}-based {N}eural {O}perators for {L}earning the {R}esponse of {C}omplex {D}ynamical {S}ystems to {L}ength-{V}ariant {M}ultiple {I}nput {F}unctions}, arXiv preprint arXiv:2311.16519,  (2023).

\bibitem{lee2023oppinn}
{\sc J.~Y. Lee, J.~Jang, and H.~J. Hwang}, {\em op{PINN}: {P}hysics-informed neural network with operator learning to approximate solutions to the {F}okker-{P}lanck-{L}andau equation}, Journal of Computational Physics, 480 (2023), p.~112031.

\bibitem{lee2024structure}
{\sc J.~Y. Lee, S.~Schotth{\"o}fer, T.~Xiao, S.~Krumscheid, and M.~Frank}, {\em Structure-{P}reserving {O}perator {L}earning: {M}odeling the {C}ollision {O}perator of {K}inetic {E}quations}, arXiv preprint arXiv:2402.16613,  (2024).

\bibitem{li2022model}
{\sc H.~Li, S.~Jiang, W.~Sun, L.~Xu, and G.~Zhou}, {\em A {M}odel-{D}ata {A}symptotic-{P}reserving {N}eural {N}etwork {M}ethod {B}ased on {M}icro-{M}acro {D}ecomposition for {G}ray {R}adiative {T}ransfer {E}quations}, Communications in Computational Physics, 35 (2024), pp.~1155--1193.

\bibitem{li2022apfos}
{\sc L.~Li and C.~Yang}, {\em {APFOS}-{N}et: Asymptotic preserving scheme for anisotropic elliptic equations with deep neural network}, Journal of Computational Physics, 453 (2022), p.~110958.

\bibitem{li2020multipole}
{\sc Z.~Li, N.~Kovachki, K.~Azizzadenesheli, B.~Liu, A.~Stuart, K.~Bhattacharya, and A.~Anandkumar}, {\em Multipole graph neural operator for parametric partial differential equations}, Advances in Neural Information Processing Systems, 33 (2020), pp.~6755--6766.

\bibitem{li2020fourier}
{\sc Z.~Li, N.~B. Kovachki, K.~Azizzadenesheli, B.~liu, K.~Bhattacharya, A.~Stuart, and A.~Anandkumar}, {\em Fourier {N}eural {O}perator for {P}arametric {P}artial {D}ifferential {E}quations}, in International Conference on Learning Representations, 2021.

\bibitem{li2024physics}
{\sc Z.~Li, H.~Zheng, N.~Kovachki, D.~Jin, H.~Chen, B.~Liu, K.~Azizzadenesheli, and A.~Anandkumar}, {\em Physics-informed neural operator for learning partial differential equations}, ACM/JMS Journal of Data Science, 1 (2024), pp.~1--27.

\bibitem{liboff2003kinetic}
{\sc R.~L. Liboff}, {\em Kinetic theory: classical, quantum, and relativistic descriptions}, Springer Science \& Business Media, 2003.

\bibitem{liu2020multi}
{\sc Z.~Liu, C.~Wei, and X.~Z.-Q. John}, {\em {M}ulti-scale {D}eep {N}eural {N}etwork ({M}scale{DNN}) for {S}olving {P}oisson-{B}oltzmann {E}quation in {C}omplex {D}omains}, Communications in Computational Physics, 28 (2020), pp.~1970--2001.

\bibitem{livadiotis2014electrostatic}
{\sc G.~Livadiotis and D.~J. McComas}, {\em Electrostatic shielding in plasmas and the physical meaning of the {D}ebye length}, Journal of Plasma Physics, 80 (2014), pp.~341--378.

\bibitem{lou2021physics}
{\sc Q.~Lou, X.~Meng, and G.~E. Karniadakis}, {\em Physics-informed neural networks for solving forward and inverse flow problems via the {B}oltzmann-{BGK} formulation}, Journal of Computational Physics, 447 (2021), p.~110676.

\bibitem{lu2021learning}
{\sc L.~Lu, P.~Jin, G.~Pang, Z.~Zhang, and G.~E. Karniadakis}, {\em Learning nonlinear operators via {D}eep{ON}et based on the universal approximation theorem of operators}, Nature machine intelligence, 3 (2021), pp.~218--229.

\bibitem{lu2022comprehensive}
{\sc L.~Lu, X.~Meng, S.~Cai, Z.~Mao, S.~Goswami, Z.~Zhang, and G.~E. Karniadakis}, {\em A comprehensive and fair comparison of two neural operators (with practical extensions) based on fair data}, Computer Methods in Applied Mechanics and Engineering, 393 (2022), p.~114778.

\bibitem{lu2022solving}
{\sc Y.~Lu, L.~Wang, and W.~Xu}, {\em Solving multiscale steady radiative transfer equation using neural networks with uniform stability}, Research in the Mathematical Sciences, 9 (2022), p.~45.

\bibitem{markowich2012semiconductor}
{\sc P.~A. Markowich, C.~A. Ringhofer, and C.~Schmeiser}, {\em Semiconductor equations}, Springer Science \& Business Media, 2012.

\bibitem{mouhot2011landau}
{\sc C.~Mouhot and C.~Villani}, {\em On {L}andau damping}, Acta Mathematica, 207 (2011), pp.~29--201.

\bibitem{nieto2001high}
{\sc J.~Nieto, F.~Poupaud, and J.~Soler}, {\em High-field limit for the {V}lasov-{P}oisson-{F}okker-{P}lanck system}, Archive for rational mechanics and analysis, 158 (2001), pp.~29--59.

\bibitem{poupaud1991diffusion}
{\sc F.~Poupaud}, {\em Diffusion approximation of the linear semiconductor {B}oltzmann equation: analysis of boundary layers}, Asymptotic analysis, 4 (1991), pp.~293--317.

\bibitem{poupaud1992runaway}
{\sc F.~Poupaud}, {\em Runaway phenomena and fluid approximation under high fields in semiconductor kinetic theory}, ZAMM-Journal of Applied Mathematics and Mechanics/Zeitschrift f{\"u}r Angewandte Mathematik und Mechanik, 72 (1992), pp.~359--372.

\bibitem{prechelt2002early}
{\sc L.~Prechelt}, {\em Early stopping-but when?}, in Neural Networks: Tricks of the trade, Springer, 2002, pp.~55--69.

\bibitem{rahaman2019spectral}
{\sc N.~Rahaman, A.~Baratin, D.~Arpit, F.~Draxler, M.~Lin, F.~Hamprecht, Y.~Bengio, and A.~Courville}, {\em On the spectral bias of neural networks}, in International Conference on Machine Learning, PMLR, 2019, pp.~5301--5310.

\bibitem{ramachandran2017searching}
{\sc P.~Ramachandran, B.~Zoph, and Q.~V. Le}, {\em Searching for activation functions}, arXiv preprint arXiv:1710.05941,  (2017).

\bibitem{stix1962theory}
{\sc T.~H. Stix}, {\em Waves in plasmas}, Springer Science \& Business Media, 1992.

\bibitem{stratton1957influence}
{\sc R.~Stratton}, {\em The influence of interelectronic collisions on conduction and breakdown in covalent semi-conductors}, Proceedings of the Royal Society of London. Series A. Mathematical and Physical Sciences, 242 (1957), pp.~355--373.

\bibitem{stratton1962diffusion}
\leavevmode\vrule height 2pt depth -1.6pt width 23pt, {\em Diffusion of hot and cold electrons in semiconductor barriers}, Physical Review, 126 (1962), p.~2002.

\bibitem{wang2023long}
{\sc S.~Wang and P.~Perdikaris}, {\em Long-time integration of parametric evolution equations with physics-informed {D}eep{ON}ets}, Journal of Computational Physics, 475 (2023), p.~111855.

\bibitem{wang2021understanding}
{\sc S.~Wang, Y.~Teng, and P.~Perdikaris}, {\em Understanding and mitigating gradient flow pathologies in physics-informed neural networks}, SIAM Journal on Scientific Computing, 43 (2021), pp.~A3055--A3081.

\bibitem{wang2021learning}
{\sc S.~Wang, H.~Wang, and P.~Perdikaris}, {\em Learning the solution operator of parametric partial differential equations with physics-informed {D}eep{ON}ets}, Science advances, 7 (2021), p.~eabi8605.

\bibitem{wang2024spi}
{\sc X.~Wang, P.~Li, K.~Jia, S.~Zhang, C.~Li, B.~Wu, Y.~Dong, and D.~Lu}, {\em {SPI}-{MION}et for surrogate modeling in phase-field hydraulic fracturing}, Computer Methods in Applied Mechanics and Engineering, 427 (2024), p.~117054.

\bibitem{wu2024capturing}
{\sc K.~Wu, X.-B. Yan, S.~Jin, and Z.~Ma}, {\em Capturing the diffusive behavior of the multiscale linear transport equations by {A}symptotic-{P}reserving {C}onvolutional {D}eep{ON}ets}, Computer Methods in Applied Mechanics and Engineering, 418 (2024), p.~116531.

\bibitem{xiao2023relaxnet}
{\sc T.~Xiao and M.~Frank}, {\em {R}elax{N}et: A structure-preserving neural network to approximate the {B}oltzmann collision operator}, Journal of Computational Physics, 490 (2023), p.~112317.

\bibitem{xu2023transfer}
{\sc W.~Xu, Y.~Lu, and L.~Wang}, {\em Transfer {L}earning {E}nhanced {D}eep{ON}et for {L}ong-{T}ime {P}rediction of {E}volution {E}quations}, Proceedings of the AAAI Conference on Artificial Intelligence, 37 (2023), pp.~10629--10636.

\bibitem{yao2007early}
{\sc Y.~Yao, L.~Rosasco, and A.~Caponnetto}, {\em On early stopping in gradient descent learning}, Constructive Approximation, 26 (2007), pp.~289--315.

\bibitem{you2022nonlocal}
{\sc H.~You, Y.~Yu, M.~D'Elia, T.~Gao, and S.~Silling}, {\em Nonlocal kernel network ({NKN}): A stable and resolution-independent deep neural network}, Journal of Computational Physics, 469 (2022), p.~111536.

\bibitem{zheng2023state}
{\sc Q.~Zheng, X.~Yin, and D.~Zhang}, {\em State-space modeling for electrochemical performance of {L}i-ion batteries with physics-informed deep operator networks}, Journal of Energy Storage, 73 (2023), p.~109244.

\end{thebibliography}

\newpage 
\appendix
\renewcommand{\thefigure}{\Alph{figure}}
\setcounter{figure}{0}

\renewcommand{\thetable}{\Alph{table}}
\setcounter{table}{0}

\section{Training History}
\begin{figure}[!htb]
	\centering 
	\subfigbottomskip=5pt 
	\subfigcapskip=-5pt 
	\subfigure[Landau daming]{
	\includegraphics[width=0.47\linewidth]{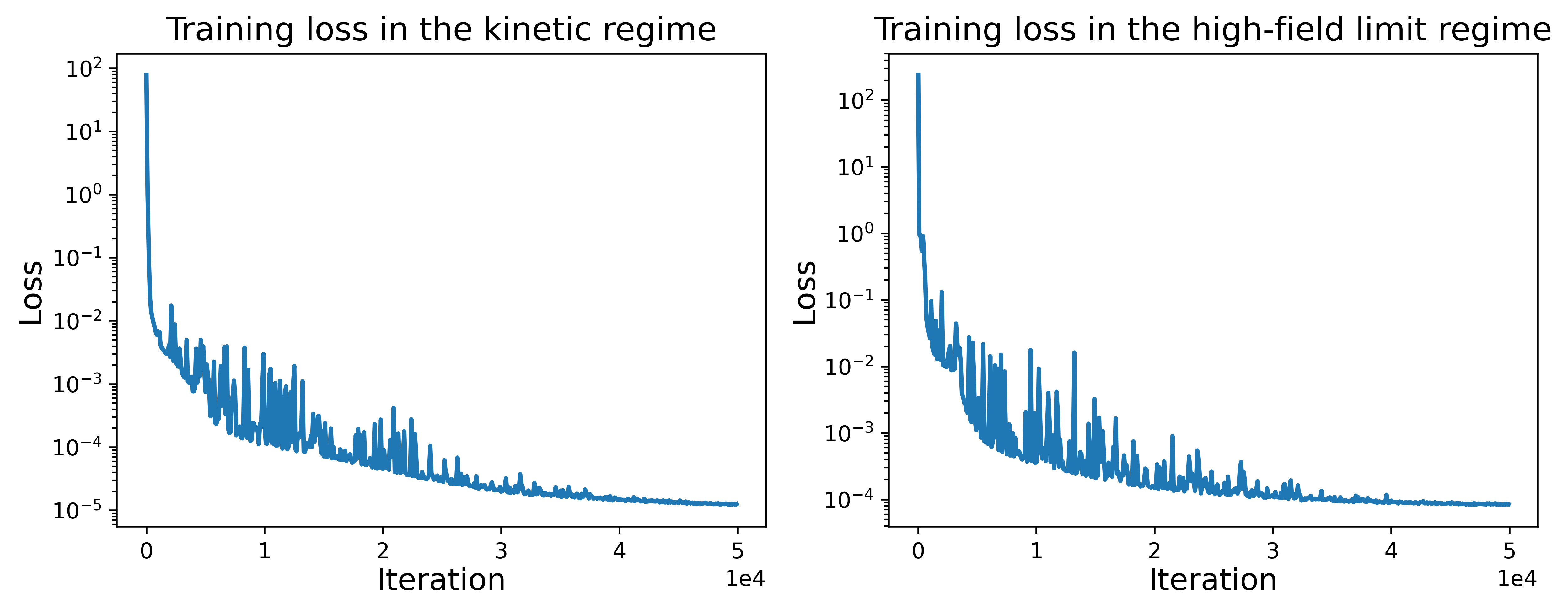}}
    \subfigure[Double peak instability]{
	\includegraphics[width=0.47\linewidth]{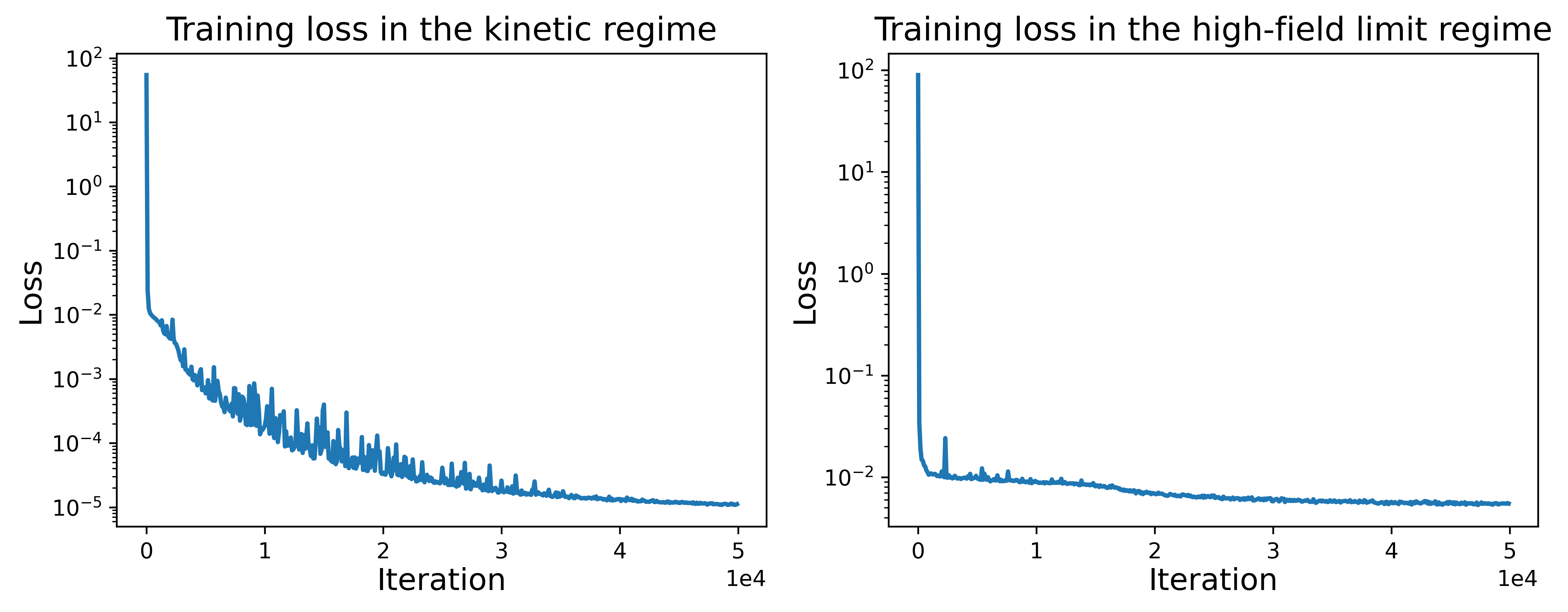}}

    \subfigure[Two stream instability]{
	\includegraphics[width=0.47\linewidth]{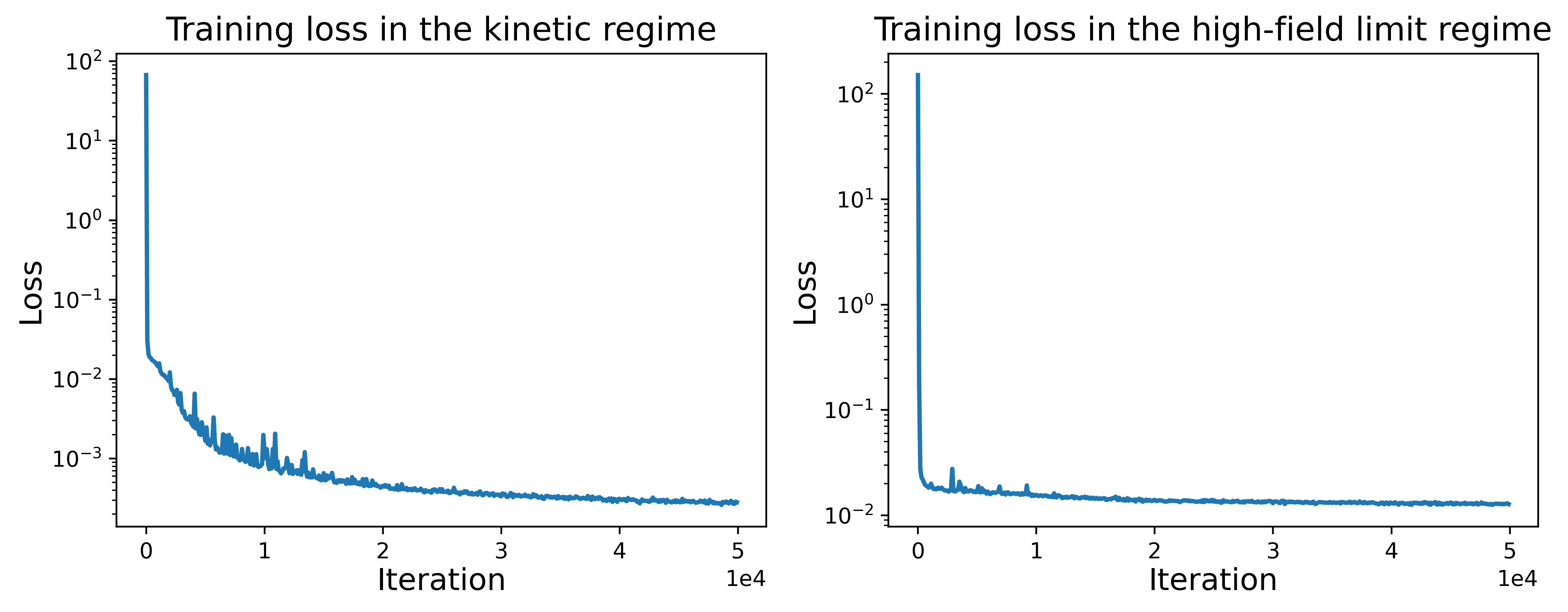}}
    \subfigure[Bump-on-tail instability]{
	\includegraphics[width=0.47\linewidth]{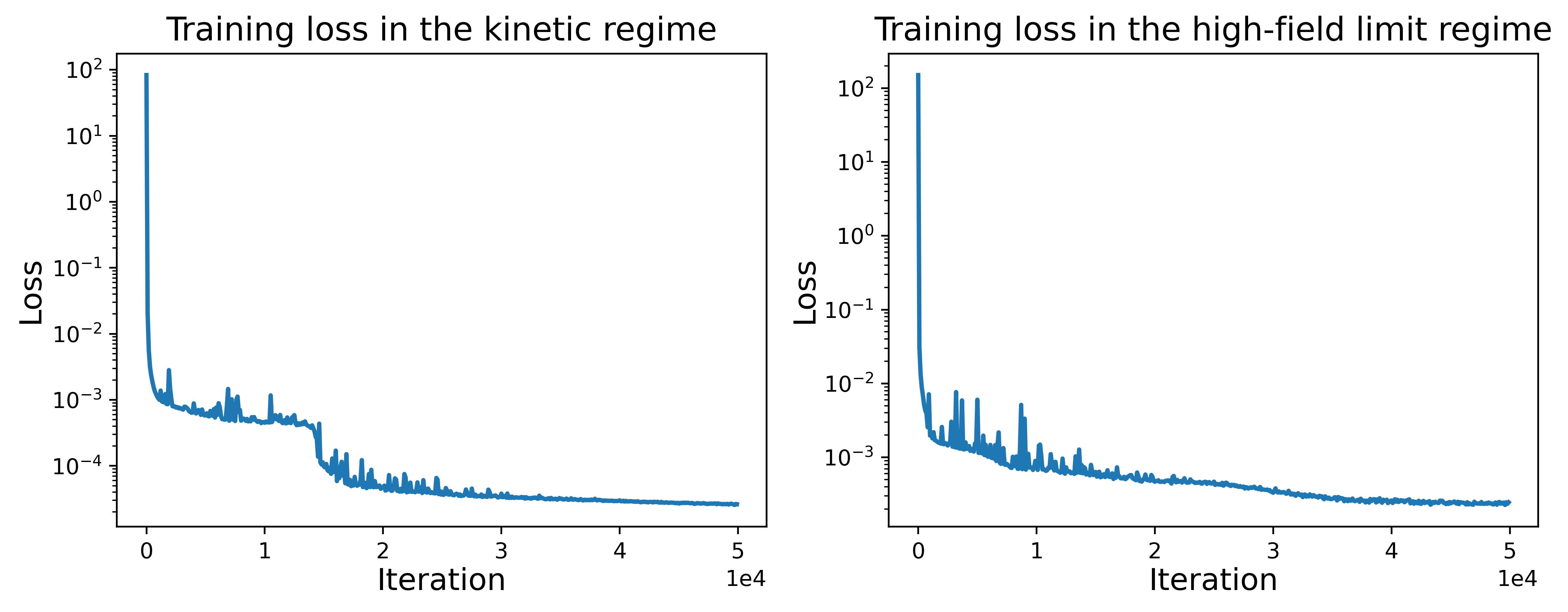}}

    \subfigure[Mixing regimes problem]{
	\includegraphics[width=0.24\linewidth]{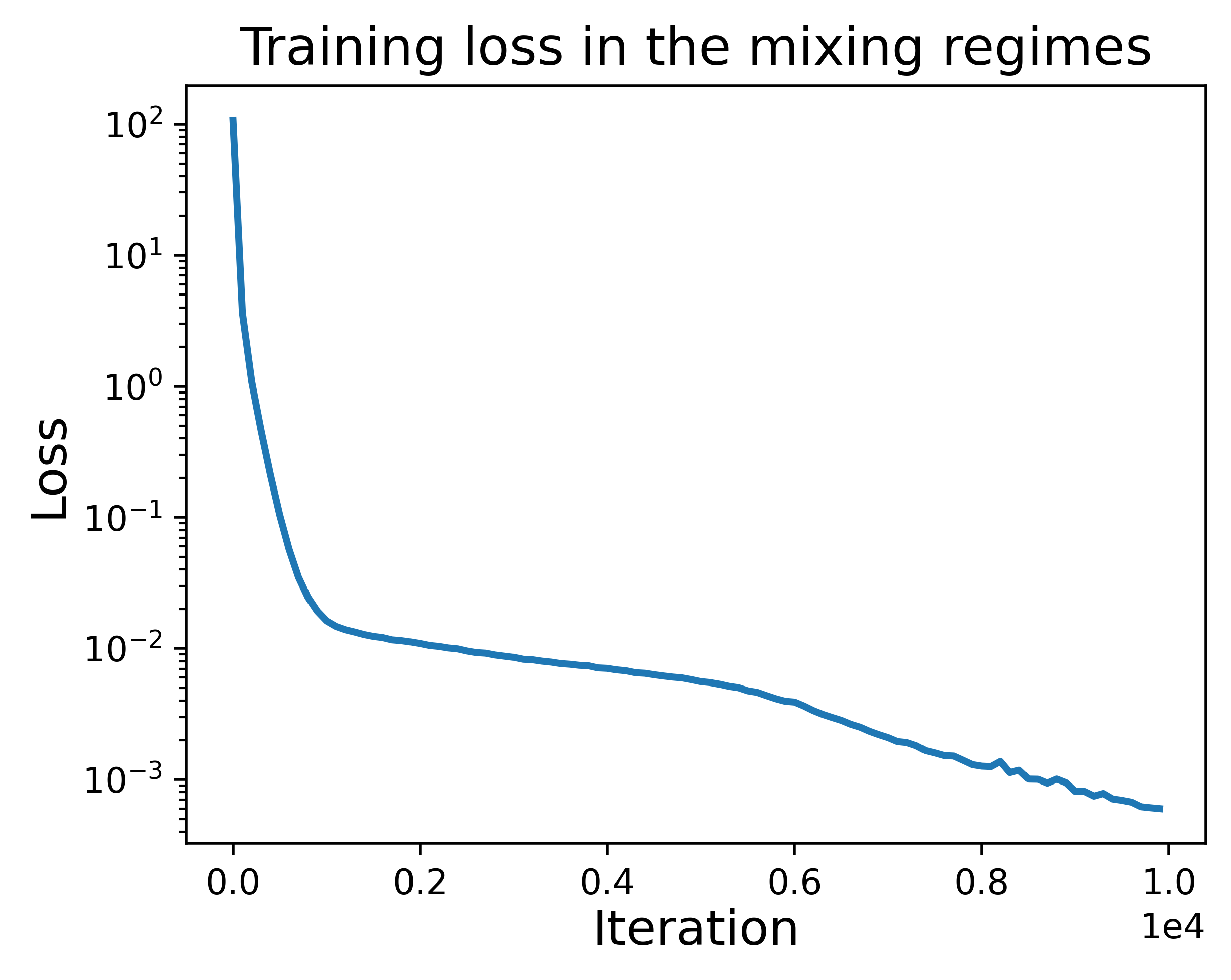}}

	\caption{Training histories of the AP-MIONet for the numerical experiments in Sect. 5. (a) Training loss convergence of the Landau damping for the VPFP system in the kinetic regime (left) and the high-field limit regime (right). (b) Training loss convergence of the double peak instability for the non-degenerate isotropic semiconductor Boltzmann equation in the kinetic regime (left) and the high-field limit regime (right). (c) Training loss convergence of the two stream instability for the non-degenerate anisotropic case in the kinetic regime (left) and the high-field limit regime (right). (d) Training loss convergence of the bump-on-tail instability for the degenerate case in the kinetic regime (left) and the high-field limit regime (right). (e) Training loss convergence of the mixing regimes problem for the degenerate case.}
    \label{Loss}
\end{figure}

Fig.\ref{Loss} summarizes the training losses of the AP-MIONet for all numerical experiments in Sect. 5.

\section{Computation cost for training process}
Table \ref{Cost} summarizes the computational time required for training the PI-MIONet and AP-MIONet models for the numerical experiments in Sect. 5.

\begin{table}[!htb]
    \centering
    \caption{Computational cost for training the PI-MIONet and AP-MIONet models for numerical experiments in Sect. 5. The training timings are obtained through parallel computation using four NVIDIA TITAN Xp GPUs.}
    \label{Cost}
    \begin{tabular}{lccc}
        \toprule
        Case  & Model (Architecture) & Regime  & Training time (hours) \\
        \midrule 
        Landau damping & PI-MIONet & Kinetic & 16.77 \\
        \midrule
        \multirow{2}{*}{Landau damping} & \multirow{2}{*}{AP-MIONet} & Kinetic & 17.02 \\
                                                   & & High-field & 17.75 \\
        \midrule
        \multirow{2}{*}{Double peak instability} & \multirow{2}{*}{AP-MIONet} & Kinetic & 17.15 \\
                                                 & & High-field & 17.21 \\
        \midrule
        \multirow{2}{*}{Two stream instability} & \multirow{2}{*}{AP-MIONet} & Kinetic & 17.84 \\
                                                   & & High-field & 17.31 \\
        \midrule
        \multirow{2}{*}{Bump-on-tail instability} & \multirow{2}{*}{AP-MIONet} & Kinetic & 19.02 \\
                                                   & & High-field & 18.43 \\
        \midrule
        Mixing regimes & AP-MIONet & Mixing & 3.51 \\
        \bottomrule
    \end{tabular}
\end{table}

\end{document}